\documentclass[11pt]{article}

\usepackage{amsmath,amssymb,amsfonts,amsthm} 
\usepackage[top=1in,bottom=1.4in,right=1in,left=1in]{geometry}

\usepackage[maxbibnames=99,maxcitenames=2,style=alphabetic,backend=biber,sorting=ynt,sortcites,natbib=true,giveninits=true]{biblatex}

\usepackage{savesym}
\savesymbol{NO}
\usepackage{graphicx}
\usepackage{hyperref}
\hypersetup{colorlinks=true,citecolor=blue,linkcolor=red}
\usepackage{epstopdf}
\usepackage{pdfpages}
\usepackage[noend]{algpseudocode}
\usepackage{algorithm}
\usepackage{algorithmicx}
\usepackage{array}
\usepackage{subcaption}
\usepackage{tikz}
\usetikzlibrary{arrows,positioning,automata}
\usepackage[mathscr]{euscript}
\usepackage{booktabs} 
\usepackage{multirow}
\usepackage{multicol}
\usepackage{longtable}
\usepackage{lscape}
\usepackage{color}
\usepackage{xifthen}
\usepackage{ifthen}
\usepackage{braket}
\usepackage{xargs}

\DeclareMathOperator*{\argmin}{argmin}

\newtheorem{proposition}{Proposition}

\newcommand{\bb}{BB}

\newcommand{\su}[2][1]{\ensuremath{\displaystyle{\sum_{#1}^{#2}}}}

\newcommand{\R}{\mathbb{R}}
\newcommand{\Z}{\mathbb{Z}}
\newcommand{\T}{\mathbb{T}}

\newcommand{\dom}{\succ}

\renewcommand{\xi}{f}
\newcommand{\U}{\mathcal{U}}
\renewcommand{\L}{\mathcal{L}}

\newcommand{\N}{\mathcal{N}}
\newcommand{\OS}{\mathcal{OS}}
\renewcommand{\top}{}

\newcommand{\bul}{\scalebox{.45}{\ensuremath{\bullet}}}
\newcommand*{\tran}{^{\mkern-1.5mu\mathsf{T}}}
\renewcommand{\top}{\tran}

\newcommand{\objf}{\boldsymbol{f}}
\newcommand{\nonneg}{\R^{2}_{\ge 0}}
\newcommand{\nd}[1]{\mathcal{ND}(#1)}

\newcommand{\revision}[1]{#1}
\newcommand{\instance}[1]{\texttt{#1}}
\newcommand{\errorone}{$\vartriangle$}
\newcommand{\errortwo}{$\blacktriangle$}
\newcommand{\percentGap}[1]{#1\%}

\renewcommand{\epsilon}{\varepsilon}

\AtEveryBibitem{
	\ifentrytype{article}{\clearfield{url}}{}
	\ifentrytype{book}{\clearfield{doi}}{}
    \iffieldundef{pages}{}{\clearfield{doi}\clearfield{eprint}}
    \ifentrytype{unpublished}{\iflistundef{publisher}{}{\clearfield{url}}}{\clearfield{url}}
    \clearfield{isbn}
    \clearfield{issn}
}

\newcommand{\keywords}[1]{\medskip\\
\noindent\textbf{Keywords.} #1}
\newcommand{\subclass}[1]{\smallskip\\
\noindent\textbf{AMS 2020 subject classification.} #1}

\title{Branch-and-bound for biobjective mixed-integer linear programming\thanks{This work was initiated when the two authors were in the School of Mathematical and Statistical Sciences at Clemson University, USA, and they were partially supported by ONR grant N00014-16-1-2725.} \textsuperscript{ }\thanks{This work used the Extreme Science and Engineering Discovery Environment (XSEDE), which is supported by National Science Foundation grant number ACI-1548562.}}

\author{Nathan Adelgren\thanks{Department of Mathematics and Computer Science, Edinboro University, PA, USA. \emph{Email}: \texttt{nadelgren@edinboro.edu}} \and Akshay Gupte\thanks{School of Mathematics, University of Edinburgh, Edinburgh, UK. \emph{Email}: \texttt{akshay.gupte@ed.ac.uk}}}

\begin{document}

\maketitle
{
\newcommand{\sep}{$\cdot$ }
\begin{abstract}
We present a generic branch-and-bound algorithm for finding all the Pareto solutions of a biobjective mixed-integer linear program. The main contributions are new algorithms for obtaining dual bounds at a node, for checking node fathoming, presolve and duality gap measurement. Our branch-and-bound is \revision{predominantly} a decision space search method since the branching is performed on the decision variables, akin to single objective problems, \revision{although we also sometimes split gaps and branch in the objective space}. The various algorithms are implemented using a data structure for storing Pareto sets. Computational experiments are carried out on literature instances \revision{and also on a new set of instances that we generate using the MIPLIB benchmark library for single objective problems}. We also perform comparisons against the triangle splitting method from literature, which is an objective space search algorithm.
\keywords{Branch-and-bound \sep Mixed-integer \sep Multiobjective optimization \sep Pareto solutions \sep Fathoming rules}
\subclass{90C11 \sep 90C57 \sep 90C29 \sep 90-04 \sep 90-08}
\end{abstract}
}

\section{Introduction}

We present a branch-and-bound (BB) algorithm that computes the \revision{nondominated solutions} of a biobjective mixed-integer linear program (BOMILP), formulated as
\begin{equation}\label{BOMILP}
\min_{x} \left\{\begin{array}{c} f_1(x) := {c^{1}}^{\top}x \\ f_2(x) := {c^{2}}^{\top}x \end{array}\right\} \quad \text{s.t.} \quad  x \in X_I := \left\{x\in \Z^{n}_+ \times \R^{p}_{+} \colon Ax\leq b,\; l_i \leq x_i \leq u_i \ \forall i\right\}.
\end{equation}
The only assumption we make on the above model is a mild and standard one: that $X_{I}\neq\emptyset$ and $-\infty < l_{i} < u_{i} < +\infty$ for all $i$, in order to have a bounded feasible problem.

BOMILPs belong to the general class of multiobjective optimization \cite{ehrgott2005multicriteria} and are an extension of the single objective mixed-integer linear program (MILP) that has been studied for decades.  
A multiobjective problem is considered solved when the entire set of \revision{nondominated} solutions has been discovered. A common approach to find these nondominated points has been to scalarize the vector objective (cf. \cite{ehrgott2006discussion,burachik2017new}) either by aggregating all objectives into one or by moving all but one objective to the constraints, but doing so does not generate all the nondominated points and supplies a very small part of the optimality information that can otherwise be supplied by the original multiobjective problem. Indeed, it is easy to construct examples of biobjective MILPs where many  Pareto solutions are located in the interior of the convex hull of the feasible set, a phenomenon that is impossible with optimal solutions of MILPs. The set of Pareto solutions of a mixed-integer multiobjective problem with a bounded feasible region is equal to the union of the set of Pareto solutions from each slice problem. Here the union is taken over the set of integer feasible values and a slice problem is a continuous multiobjective program obtained by fixing the integer variables to some feasible values. In general, there could be exponentially \revision{many slice problems, each providing some} Pareto solutions. Enumeration of the Pareto set for a pure integer problem has received considerable attention, including iterative approaches  \cite{lokman,ozpeynirci2010exact} and lower and upper bounds on the number of Pareto solutions \cite{bazgan2013number,stanojevic2013cardinality} under certain assumptions. Algorithms using rational generating functions to enumerate all the nondominated optima in polynomial-time for fixed parameters (either size of decision space or number of objectives) were given in \cite{de2009pareto,blanco2012complexity}.
There also have been many efforts at finding good approximations of the Pareto set \cite{bazgan2015approximate,grandoni2014new,ruzika2005approximation,sayin2000measuring,sayin2003procedure,mittal2013general}.

\subsection{Background on existing methods} 

Algorithms for exact solution of multiobjective mixed-integer problems (MOMILPs) can be broadly classified into three categories depending on the underlying techniques they use: (i) those based on scalarization methods that transform the MOMILP into a MILP with a modified objective or with new constraints, (ii) branch-and-bound algorithms which are decision space search since they divide the feasible region $X_{I}$ by branching on variables (in a manner similar to solving MILPs), and (iii) those based on objective/criterion space search methods that solve MILPs or multiobjective LPs over subsets of the feasible objective space $f(X_{I}) := \{(f_{1}(x), f_{2}(x)) \colon x \in X_{I} \}$. Multiobjective pure integer problems have been extensively studied in literature and many scalarization methods have been developed, either specifically for biobjective problems \cite{ralphs2006improved} or the fully general multiobjective case \cite{ehrgott2008improved,mavrotas2013improved,przybylski2010149,ehrgott2006discussion}. Specific classes of biobjective combinatorial problems also have algorithms for solving them \cite{sourd2008multiobjective,leitner2014computational,royset2007solving,biobjkp1998,berube2009exact}. Earliest branch-and-bound methods for multiobjective pure integer programs can be found in  \cite{kiziltan1983algorithm,KLEIN1982378}, but since then more sophisticated algorithms have been developed \cite{parragh2015boip,jozefowiez2012generic,gadegaard2019bi,belotti2012biobjective,vincent2013biobjective,stidsen2014branch}. \revision{Objective space search methods have been the focus of study in the last few years separately for pure integer programs \cite{kirlik2014new,boland2015bcriterion,boland2016quadrant,boland2015shape,dachert2015linear,turgut2019exact} and mixed-integer programs \cite{perini2020criterion,boland2015acriterion}.}

Algorithms specialized for the pure integer case do not extend to the mixed-integer case primarily because of the way they certify Pareto optimality. The Pareto set of a mixed-integer problem is a finite union of graphs of piecewise linear functions, whereas that for a pure integer problem is a finite set of points, and hence Pareto computation and certification of Pareto optimality of a given subset is far more  complicated in the former case.  In fact, mixed-integer problems can benefit immensely from sophisticated data structures for storing Pareto sets, as shown recently by \citet{treestructure}. Most of the \bb{} algorithms in literature are designed specifically for problems where all the integer variables are binary; see the literature reviews in \cite{przybylski2017multi,gadegaard2019bi}. Correct node fathoming rules are necessary to guarantee correctness of a \bb{} algorithm. \citet{belotti2015fathoming,belotti2012biobjective} have proposed sophisticated algorithms, based on solving LPs, for node fathoming rules and checking Pareto optimality, and report some limited preliminary computational results. In principle, this leads to a \bb{} algorithm for BOMILP with general integer variables, however, such an algorithm based on sophisticated node fathoming rules has neither been fully implemented nor extensively tested.

\revision{The original preprint of this paper was first made available as \cite{adelgren2016bb}. Since then, several other papers on exact methods for generating the entire nondominated set have appeared in literature \cite{boland2019preprocessing,fattahi2018one,rasmi2019gondef,soylu2016exact,soylu2018search}, none of which are based on a BB method and all of which seem to be unaware of our work. Some of these methods also have the restriction that they can handle only binary variables, not general integers. Of particular note in this regard is the recent paper of \citet{gadegaard2019bi} which although is based on a BB algorithm, does not acknowledge the algorithms already developed in \cite{adelgren2016bb,belotti2012biobjective}.
}

\subsection{Summary of our work}
Our exact algorithm for general BOMILP is based on the \bb{} method. Although there is certainly merit in studying and developing objective space search methods for solving BOMILP, our choice is motivated by the recognition that there is still much work that can be done to exploit the structure of Pareto points in biobjective problems to improve \bb{} techniques for BOMILP. That is indeed the main contribution of this paper --- an exhaustive computational study of ideas that specifically address the biobjective nature of problem~\eqref{BOMILP}. Besides the fact that \bb{} operates mainly in the $x$-space and objective space search, as the name suggests, operates solely in the $f$-space, another point of distinction between the two is that the MILPs we consider at each node of the \bb{} tree do not have to be solved to optimality whereas the correctness of the latter depends on MILPs being solved to optimality. Of course, it is to be expected that solving MILPs for a longer time will lead to better convergence results for our \bb{}. Implementing our \bb{} through the callback interface of a MILP solver allows us to utilize the huge computational progress made in different components of \bb{} for MILP (cf. \cite{achterberg2013mixed,morrison2016branch}). 

The main components of any \bb{} for MILP include presolve, preprocessing, primal heuristics, dual bounding via cutting planes, node processing, and branching. We present new algorithms to adapt and extend each of these components to the biobjective case. We begin with presolve; since primal presolve techniques work solely on the feasible region, their implementations in state-of-the-art MILP solvers can be directly used for a BOMILP. However, dual presolve utilizes information from the objective function and hence cannot be used directly for a BOMILP. We are the first to discuss (\textsection\ref{sec:presolve}) and implement an extension of a variety of dual presolve techniques to the multiobjective setting. Additionally, we show that using one of the primal presolve techniques --- probing on integer variables (\textsection\ref{sec:probe}), alongside branching reduces the overall computational time. Two different preprocessing algorithms (\textsection\ref{sec:preprocessing}) are proposed for generating good primal bounds. Our main body of work is in developing new node processing techniques (\textsection\ref{sec:np}) for BOMILP.  The node processing component takes increased importance for BOMILP since bound sets for a multiobjective problem are much more complicated than those for a single objective problem (cf. \textsection\ref{sec:bounds}), meaning that generation of valid dual bounds and fathoming of a node is not as straightforward as that for MILPs. At each node, we describe procedures to generate valid dual bounds while accounting for the challenges of biobjective problems and strengthen these bounds through the use of locally valid cutting planes and the solution of single objective MILPs. Our bounds are tighter than what has previously been proposed. To guarantee correctness of our \bb{}, we develop new fathoming rules and delineate their difference to the work of \cite{belotti2012biobjective,belotti2015fathoming} in \textsection\ref{sec:compare}. \revision{We use the depth first rule for node selection.} A branching scheme is presented in \textsection\ref{sec:branch} and a method for exploiting distances between Pareto points in the objective space is discussed in \textsection\ref{sec:osgap}. Finally, our \bb{} also incorporates an early termination feature that allows it to terminate after a prescribed gap has been attained. In the MILP case, gap computation is trivial to implement because primal and dual bounds for MILPs are scalars. However for BOMILPs, since these bounds are subsets of $\R^{2}$ as explained in \textsection\ref{sec:bounds}, computation of optimality gap requires the use of error measures that are nontrivial to compute. To aid quicker computation, we propose in \textsection\ref{sec:gap} an approximated version of the Hausdorff metric and computationally compare it to the hypervolume gap measure from literature.


An extensive computational analysis is carried out in \textsection\ref{sec:compute} 
on literature instances. 
The first of these experiments evaluates our three dual presolve techniques 
and the results show that duality fixing is the most useful of the three for reducing CPU time. In our second experiment, we demonstrate that preprocessing methods utilizing $\epsilon$-constraint scalarization techniques typically yield better primal bounds at the start of \bb{} than  weighted sum scalarization techniques. Next, we evaluate the performance of various procedures, such as probing, objective-space fathoming, a variety of cut generation techniques, and some minor improvements to our proposed fathoming rules, that we propose in this paper for improving the overall performance of \bb{}. These tests indicated that probing prior to each branching decision and objective space fathoming are very useful for decreasing the total solution time. 
The local cuts that we added were not as useful. Finally, we compared the performance of our \bb{} with that of the triangle splitting method \cite{boland2015acriterion}, which we recall is an objective space search method, and observe that our \bb{} uses less CPU time 
to compute the complete Pareto sets of the test instances. 

We conclude this paper with a few remarks in \textsection\ref{sec:conclude}. 
We observe that a majority of the algorithms proposed in this paper can be extended naturally to the multiobjective case. The main challenge in developing a fully implementable and efficient \bb{} algorithm for multiobjective MILP is in carrying out the bound domination step. We present some directions for future research on this topic.

\revision{Comparing the contributions of this paper to the literature review on other \bb{} algorithms found in \citep[Table 1]{gadegaard2019bi}, we see that there are several distinguishing features of our algorithm: 
\begin{enumerate}
\item ours is the \emph{only} full implementation so far that handles mixed-integer problems, whereas all others are either pure binary or mixed-binary;
\item branching is not simply a matter of variable fixing for us and we develop a scoring scheme for branching motivated by reliability branching in MILP; 
\item we are the only \bb{} that does extensive presolve and preprocessing at the root node; 
\item we perform probing before branching, which provides a huge speedup; 
\item for our lower bound, we employ all of the strategies, except hypersurface, used by others but also go a step further. Instead of stopping by using supported points of the BOMILP, we use dual information from the MILP to use a tighter lower bound than just the supported points themselves -- specifically, we use a line segment bound generated from the MILP's dual information;
\item we propose a new optimality gap measure, that is an approximate Hausdorff metric in the objective space, to allow for early termination of our \bb{} once a prescribed gap has been attained.
\end{enumerate}
There are also similarities with others, such as Pareto branching and use of cutting planes in the objective space.
}

\section{Preliminaries}\label{sec:prelim}


The idea of optimality for single objective optimization is replaced with the idea of \emph{efficiency} in multiobjective problems. Consider BOMILP \eqref{BOMILP}. For any two points $y,y^{\prime}\in\R^{2}$, it is said that $y$ \emph{dominates} $y^{\prime}$ if $y \le y^{\prime}$, or equivalently $y^{\prime}\in y + \nonneg$. We express this relationship as $y \dom y^{\prime}$. 
Denoting $\objf(x) := (f_{1}(x),f_{2}(x))$, which is a vector in $\R^{2}$, a point $x\in X_I$ is said to be \emph{efficient} if there is no $x^{\prime} \in X_I$ such that $\objf(x^{\prime}) \dom \objf(x)$. A point in $\R^{2}$ is Pareto optimal (also called nondominated) if it is the $\objf$-image of some efficient solution in $X_{I}$. Denote the sets of efficient solutions, \revision{also called Pareto optimal solutions}, and nondominated solutions respectively, by
\[
X_{E} := \{x\in X_{I} \colon x \text{ is efficient} \}, \qquad Y_{N} = \objf(X_{E}) := \{\objf(x) \colon x \in X_{E}\} .
\]
The nondominated subset of any $S\subset\R^{2}$ is defined as 
\[
\nd{S} := \{y\in S\colon \nexists\, y^{\prime}\in S \text{ s.t. } y^{\prime}\dom y \}.
\] 
Therefore, if we let $Y_{I} := \{\objf(x) \colon x\in X_I\}$, we have that $Y_{N} = \nd{Y_{I}}$.

For $k=1,2$, let $f^{\ast}_{k} := \min\{f_{k}(x)\colon x\in X_{I}\}$ be the optimal value of objective $k$ for the single objective problem. Denote
\[
Y^{k}_I := \left\{y\in \R^{2} \colon y_{i} = f^{\ast}_{i} \ \, i\neq k,\, y_k = \min_{x\in X_I} \left\{ f_k(x) \colon f_i(x) = f^{\ast}_{i} \ \,  i\neq k\right\} \right\} \quad k=1,2.
\]
\revision{This set is a singleton and is called the lexicographic minima of the biobjective problem of minimising the two objectives over $X_{I}$.} We have $Y^{k}_{I}\subset Y_{N}$. For each of $X_I$, $Y_I$, and $Y^k_I$, dropping the $I$ subscript indicates the continuous relaxation of the set. Also, if we add a subscript $s$, then it means that the set is associated with node $s$ of the BB tree. We use $\OS$ to denote the \emph{objective space}\footnote{\revision{Note that there is a slight deviation from literature where it is common to denote $\R^{2}$ as the objective space.}}, i.e., the smallest rectangle in $\R^2$ that contains $Y$. Given $S\ \subseteq \OS \subseteq \R^2$, the \emph{ideal point} of $S$, denoted $S^{ideal}$, is the point $y \in \R^2$ with $y_k = \min_{y \in S}\{y_k\}$ for $k =1,2$. 

We assume background in branch-and-cut algorithms for single objective problems (cf. \cite{morrison2016branch}). One of the key differences and challenging aspects of BOMILP versus MILP is the concept of primal and dual bound sets, which we explain next.

{
\newcommand{\Pareto}[1]{\mathcal{P}(#1)}
\subsection{Bound sets for BOMILP}\label{sec:bounds}
Similar to the single objective case, correct fathoming rules are essential for any \bb{} algorithm to solve BOMILP to Pareto optimality (\revision{generate the entire nondominated frontier}). Primal and dual bounds in a single objective \bb{} are scalars, making it easy to compare them and fathom a node by bound dominance. In biobjective BB, these bounds are subsets of $\R^{2}$. Bound sets were first discussed by \citet{ehrgott2007bound}. The manner in which these bound sets are generated within a \bb{} is conceptually similar to the single objective case and we explain this next. Note that our forthcoming explanation trivially extends to the multiobjective case. 

Suppose that we are currently at node $s$ of the \bb{} tree. The primal bound sets are constructed from the set of integer feasible solutions, denoted by $T_{s}\subset\Z^{n}_{+}$, found so far by the \bb{}. 
For every $\tilde{x}\in T_{s}$, the BOLP obtained by fixing $x_{i}=\tilde{x}_{i}$ for $i=1,\dots,n$ in BOMILP \eqref{BOMILP} is called the \emph{slice problem}. The Pareto curve for this slice problem is $\nd{\objf(X(\tilde{x}))}$, where $X(\tilde{x})$ denotes the feasible set of the slice problem, and this curve is convex (because it is minimization) and piecewise linear. Then $\N_{s} := \displaystyle\nd{\cup_{\tilde{x}\in T_{s}}\nd{\objf(X(\tilde{x}))}}$ is the globally valid primal bound calculated at node $s$. For the dual bound set, we consider BOLPs obtained by relaxing integrality on variables. Since $X_{s}$ denotes the relaxed feasible set at node $s$ and $Y_{s}=\objf(X_{s})$, the local dual bound is $\nd{Y_{s}}$ and is convex piecewise linear. The global dual bound $\L^{global}_{s}$ is obtained by considering the local dual bounds for all the open nodes in the \bb{} tree, i.e., $\L^{global}_{s} = \nd{\cup_{s^{\prime}\in\Omega_{s}}\L_{s^{\prime}}}$ where $\Omega_{s}$ is the set of unexplored nodes so far, and this bound is a union of convex piecewise linear curves.

For multiobjective \bb{}, node $s$ is allowed to be fathomed by bound dominance if and only if $\L_{s}$ is dominated by $\N_{s}$, i.e., for every $y^{\prime}\in\L_{s}$ there exists a $y\in\N_{s}$ such that $y \dom y^{\prime}$. Equivalently, due to translation invariance of $\dom$, we have that node $s$ can be fathomed by bound dominance if and only if $\L_{s}+\nonneg \subset \N_{s}+\nonneg$, \revision{where $\nonneg$ is the nonnegative orthant}. For this reason, henceforth for convenience, we consider our local dual bound to be $\L_{s}=\nd{Y_{s}}+\nonneg$ and the current primal bound to be $\U_{s}:=\N_{s}+\nonneg$. Thus the dual bound set is a polyhedron whereas the primal bound is a finite union of polyhedra. Although this deviates from the traditional view of bound sets, which defines them in the previous paragraph in terms of the boundary of these polyhedra, it is clear that there is a one-to-one correspondence between fathoming rules for the two alternate representations of bound sets.

Figure~\ref{bound_sets} illustrates the concept of bound sets. Here, $s_2$ can be fathomed because $\L_{s_{2}}\subset\U_{s}$ but we cannot say anything about fathoming node $s_1$ since $\L_{s_{1}}\nsubseteq\U_{s}$. As can be imagined from Figure~\ref{bound_sets}, fathoming is even more crucial and computationally expensive for BOMILPs since it involves checking inclusion and intersection of polyhedral sets as opposed to comparing scalar values in the MILP case. Thus, the majority of the computational effort in multiobjective BB is spent processing a node $s$ of the BB tree, in particular checking various fathoming rules.

\begin{figure}[h]
\centering
\begin{tikzpicture}
\node[anchor=south west,inner sep=0] (image) at (0,0) {\includegraphics[width=0.3\textwidth]{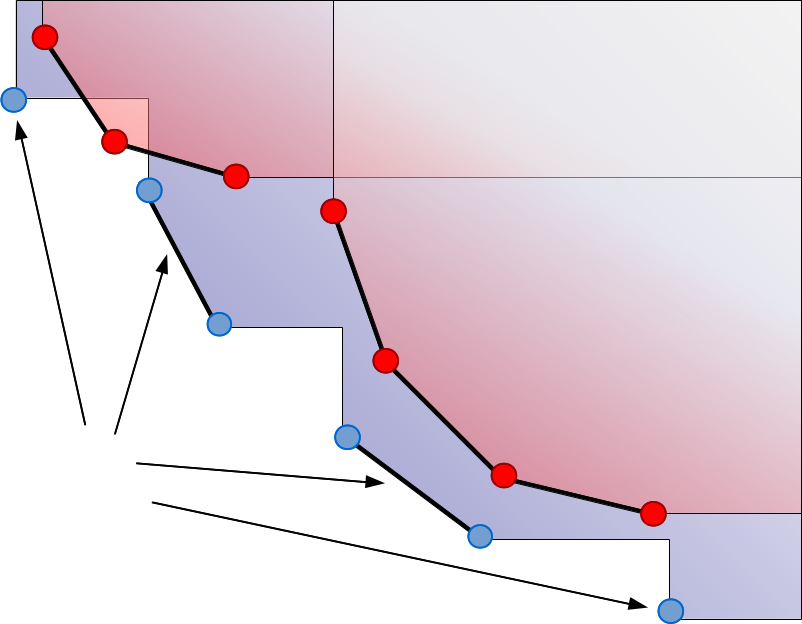}};
\node[align=center,yshift=0.35cm,xshift=-.75cm] at (image.center) {\large $\U_s$ };
\node[align=center,yshift=1.25cm,xshift=-.8cm] at (image.center) {\large $\L_{s_1}$ };
\node[align=center,yshift=0.15cm,xshift=0.95cm] at (image.center) {\large $\L_{s_2}$ };
\node[align=center,yshift=-1.05cm,xshift=-1.75cm] at (image.center) {\large $\N_s$ };
\end{tikzpicture}
\caption{Primal $(\U)$ and dual $(\L)$ bound sets for BOMILP}
\label{bound_sets}
\end{figure}
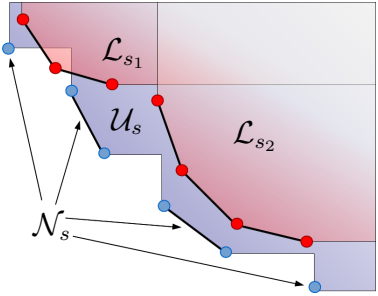

}

\section{Presolve and Preprocessing}\label{sec:pre}
Examining the structure of an instance of single objective MILP prior to solving it, and utilizing information found during this examination to simplify the structure of the instance often has had a significant impact on the time and effort needed to solve that instance. 
It has also been shown that knowledge of feasible solutions for an instance of MILP can have a significant impact on solution time. 
Hence, it seems natural as a first step to extend the techniques used in these procedures to the biobjective case. For the discussion that follows we distinguish the idea of simplifying an instance of BOMILP based on its problem structure from the idea of determining a set of initial integer feasible solutions. We refer to the first as \emph{dual presolve} and the latter as \emph{preprocessing}. 

\subsection{Dual Presolve} \label{sec:presolve}

Presolve for MILP uses both primal and dual information. The primal information of a BOMILP instance is no different than its single objective counterpart and thus primal presolve techniques can be applied directly to it. However, due to the presence of an additional objective, one must take care while utilizing dual information for a biobjective problem. 
We extend a few single objective dual presolve techniques to BOMILP (their extension to three or more  objectives is immediate and omitted here). 

First, we extend duality fixing (cf. \cite{martin2001general}). Let $a_{rj}$ denote the element of matrix $A$ in row $r$ and column $j$ and $c^{k}_{j}$ be the $j^{th}$ entry of $k^{th}$ objective. 

\begin{proposition}[Duality fixing]\label{duality_fixing}
Suppose there exists a $j$ with $c^k_j \geq 0$ and $a_{ij} \geq 0$ for all $k,i$. Then \revision{$Y_{N} \subseteq f(\{x\in X_{I}\colon x_{j} = l_{j} \})$.}
Similarly, if there exists a $j$ with $c^k_j \leq 0$ and $a_{ij} \leq 0$ for all $k,i$, then \revision{$Y_{N} \subseteq f(\{x\in X_{I}\colon x_{j} = u_{j} \})$.}
\end{proposition}
\begin{proof}
It is well known (cf. \citep[Theorem 4.5]{ehrgott2005multicriteria}) that $x^*$ is efficient for a MOMILP if and only if there exists $\epsilon$ such that $x^*$ is optimal to the problem: 
\begin{equation}\label{prop1_eqn}
\min_{x} \{f_1(x)\colon x \in X_I,\, f_k(x) \leq \epsilon_k \text{ for all } k \neq 1\}
\end{equation}
Hence, every efficient solution to the given BOMILP can be obtained by solving \eqref{prop1_eqn} for some $\epsilon$. If the stated assumptions hold, then single objective duality fixing can be applied to \eqref{prop1_eqn}. This shows that every efficient solution to the given BOMILP can be obtained by solving the modified version of \eqref{prop1_eqn} in which variable fixing has been performed.
\end{proof}

Next, we extend the exploitation of singleton and dominating columns \cite{gamrath2015}. 

\begin{proposition}[Singleton Columns]\label{singleton_columns}
\renewcommand{\gamma}{t}
For every row $r$ in the system $Ax\le b$, define ${J}(r) := \{j \colon a_{rj} > 0, c_j^k < 0 \ \forall k, a_{ij} = 0 \ \forall i\neq r\}$ and 
\[
{U}_r := \sum_{j\in {J}(r)} a_{rj}l_j + \sum_{j\not\in {J}(r), a_{rj} > 0} a_{rj}u_j + \sum_{j\not\in {J}(r), a_{rj} < 0} a_{rj}l_j.
\] 
Suppose there exists some \revision{$i\in {J}(r)$} such that \[\frac{c_i^k}{a_{ri}} \;\leq\; \min\left\{\frac{c_\gamma^k}{a_{r\gamma}} \colon \gamma \in {J}(r), \gamma \neq i\right\}.\] If $a_{ri}(u_i-l_i) \leq b_r - {U}_r$, then $X_{E}\subseteq\{x\colon x_{i}=u_{i}\}$.
\end{proposition}
\begin{proof}
\renewcommand{\ell}{l}
Let $x$ be an efficient solution with $x_i  < u_i$. If $x_j = \ell_j$ for all $j\in J (r )\setminus\{i\}$,
then a new solution $x'$ constructed from $x$ by setting $x_i'$ to $u_i$ is feasible because
\[
\su[j]{}a_{rj}x_j'=\su[j\neq i]{}a_{rj}x_j' + a_{ri}u_i \leq U_r + a_{ri}(u_i-\ell_i) \leq b_r.
\]
Additionally, the value of every objective function improves because $c_i^k < 0$ for all $k$. This contradicts
our assumption of $x$ being efficient. Hence, there exists a $j \in J (r )\setminus\{i\}$ with $x_j > \ell_j$. In this case we can construct a new solution ${x}^*$ from $x$ by decreasing the value of $x_j$
to $x_j'$ while at the same time increasing the value of $x_i$ so that $A_{r \bul} {x}^* = A_{r\bul} x$. In particular, $a_{r i} ({x}_i^* - x_i) = a_{r j} (x_j - {x}_j^* )$ holds. The change of objective $k$ can be estimated by
\begin{align*}
c_i^k {x}_i^* + c^k_j {x}_j^* &= c^k_i x_i + c^k_j x_j + c^k_i ({x}_i^* - x_i) - c^k_j (x_j - {x}_j^* )\\
&= c^k_i x_i + c^k_j x_j + c^k_i \frac{a_{rs}}{a_{rs}} ({x}_i^* - x_i) - c^k_j \frac{a_{rj}}{a_{rj}} (x_j - {x}_j^* )\\
&\leq c^k_i x_i + c^k_j x_j + c^k_i \frac{a_{rs}}{a_{rs}} ({x}_i^* - x_i) - c^k_i \frac{a_{rj}}{a_{rs}} (x_j - {x}_j^* )\\
&= c^k_i x_i + c^k_j x_j + \frac{ c^k_i}{a_{rs}}\left( a_{rs}({x}_i^* - x_i) - a_{rj} (x_j - {x}_j^* )\right)\\
&= c^k_i x_i + c^k_j x_j.
\end{align*}
If ${x}_i^* = u_i$, the result of the proposition holds. Otherwise, ${x}_j^* = \ell_j$ holds. Applying this argument iteratively results in an optimal solution with ${x}_i^* = u_i$ or ${x}_j^* = j$ for all $j \in J (r )\setminus\{s\}$. But as shown before, the latter case contradicts the efficiency of $x^*$.
\end{proof}

A similar procedure can be followed for $a_{rj} < 0$, $c_j^k > 0$ for all $k$, thereby fixing $x_{s}=l_{s}$. 

Given two variables $x_{i}$ and $x_{j}$, either both integer or both continuous, we say that $x_j$ \emph{dominates} $x_i$ if (i) $c_j^k \leq c^k_i$ for all $k$, and (ii) $a_{r j} \leq a_{ri}$ for every $r$. \footnote{This variable domination has no relationship with the idea of domination between bound sets}

\begin{proposition}[Dominating columns]\label{dominating_col_disjunction}
If $x_j$ dominates $x_i$,  
\[
Y_{N} \,=\, \left\{ \objf(x)\colon x \in X_{E},\, x_{i}=l_{i} \text{ or } x_{j}=u_{j}\right\} \, \subseteq \, \left\{ \objf(x)\colon x \in X_{I},\, x_{i}=l_{i} \text{ or } x_{j}=u_{j}\right\}.
\]
\end{proposition}
\begin{proof}
The $\subseteq$-inclusion is obvious from $X_{E}\subseteq X_{I}$, and so we have to argue the equality. We will need the following claim, which can be argued easily and is also an extension of \citep[Lemma 1]{gamrath2015} : for any $x\in X_{I}$ with a pair of indices $(i,j)$ such that $x_{j} < u_{j}$, $x_{i} > l_{i}$, and $x_j$ dominates $x_i$, the point $x^{\alpha}$ constructed for arbitrary $0 < \alpha \le \min\{ x_i - l_i, u_j - x_j \}$ as follows,
\begin{equation}		\label{claimpoint}
x^{\alpha}_{i} = x_{i} - \alpha, \qquad x^{\alpha}_{j} = x_{j} + \alpha, \qquad x^{\alpha}_{t} = x_{t}, \quad t\neq i,j,
\end{equation}
satisfies $x^{\alpha}\in X_{I}$ and $f_{k}(x^{\alpha}) \le f_{k}(x)$ for all $k$.

Since $Y_{N} = \objf(X_{E})$ by definition, the $\supseteq$-inclusion is obvious. Now suppose for sake of contradiction that the $\subseteq$-inclusion is not true. Then there exists some $y\in Y_{N}$ for which 
\begin{equation}		\label{contradict}
\objf^{-1}(y) \, \bigcap \, \left( \{x\colon x_{i}=l_{i} \} \cup \{x \colon x_{j}=u_{j} \} \right) = \emptyset.
\end{equation}
Take any $x\in \objf^{-1}(y)$, this point has $x_j < u_j$ and $x_i > l_i$. Consider the feasible solution $x^\alpha$, for $\alpha = \min\{ x_i - l_i, u_j - x_j \}$, constructed as in equation~\eqref{claimpoint}. By definition of $\alpha$, we have $x^\alpha_j = u_j$ or $x_i^\alpha = l_i$, and the claim gives us $x^{\alpha}\in X_{I}$. We know that $\objf(x) = y\in Y_{N}$. Then, $\objf(x^{\alpha}) \le \objf(x)$ from the above claim  implies that $x^{\alpha}\in\objf^{-1}(y)$. Hence, we have reached a contradiction to equation~\eqref{contradict}. 
\end{proof}

One may use the disjunction resulting from Proposition \ref{dominating_col_disjunction} to generate valid cutting planes for $X_I$ prior to branching. Additionally, there are also ways to further utilize the structure of dominating columns in order to strengthen variable bounds as described in \citet[Theorem 3, Corollary 1 and 2]{gamrath2015}. These methods for strengthening bounds also extend to the multiobjective case. However, we did not find these methods to be advantageous in our experiments. Thus, since the description of these additional strategies is quite lengthy, we omit them from this work.

\subsection{Preprocessing} \label{sec:preprocessing}  

As in the single objective case, the efficiency of BB can be significantly improved if good-quality primal feasible solutions can be generated prior to the start of BB. This can be accomplished by a heuristic method, such as \citep{soylu2015heuristic,leitner2016ilp}. We utilize two different preprocessing techniques, both of which solve single objective MILPs subject to a certain time limitation --- the first uses the $\epsilon$-constraint method, and the second uses the weighted-sum approach. We briefly discuss the benefits and drawbacks of using either the $\epsilon$-constraint or weighted-sum approaches (see \citep{ehrgott2005multicriteria} for background on scalarization methods).

\paragraph{\textbf{$\epsilon$-constraint:}} It is well known that for a BOMILP every $y \in Y_N$ can be obtained using the $\epsilon$-constraint method. Unfortunately though, when a MILP formulated using the $\epsilon$-constraint method is not solved to optimality, there are two major drawbacks: (i) each $y\in Y_I$ discovered while processing the MILP must lie within a restricted region of $\OS$, and (ii) the information associated with the best dual bound cannot be utilized. 

\paragraph{\textbf{weighted-sum:}} The major drawback of the weighted sum method is that when a MILP is formulated using this method, only \emph{supported} Pareto solutions can be found, i.e., those lying on the boundary of the convex hull of $Y_N$. There are, however, the following two benefits: (i) $y \in Y_I$ discovered during the MILP solve are not restricted to any particular region of $\OS$, and (ii) the best dual bound discovered during the MILP solve is valid for all $y \in Y_I$ and can therefore be used to create a cutting plane in $\OS$.

\begin{figure}
\begin{subfigure}[h]{.48\textwidth}
\centering
\includegraphics[width=0.8\textwidth]{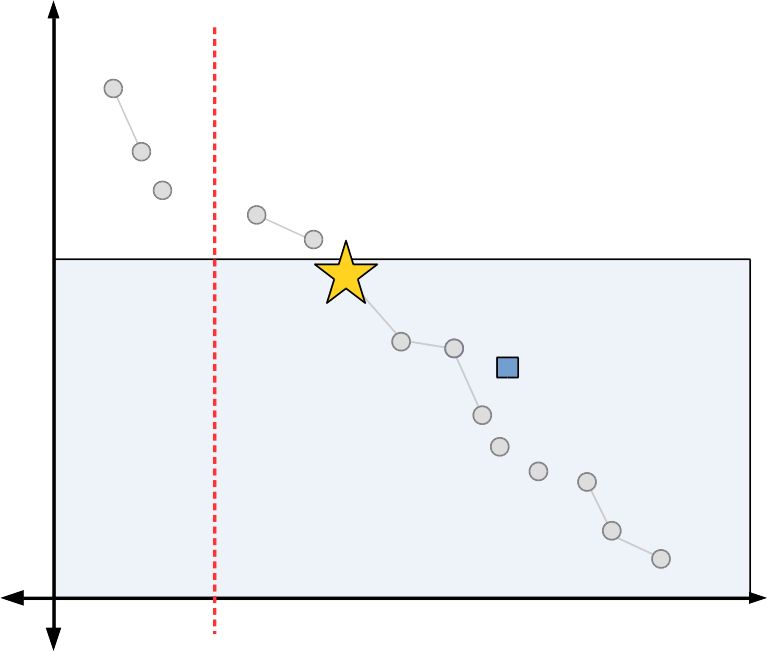}
\caption{$\epsilon$-constraint method}
\label{preprocessing_epsilon}
\end{subfigure} 
\begin{subfigure}[h]{.48\textwidth}
\centering
\includegraphics[width=.8\textwidth]{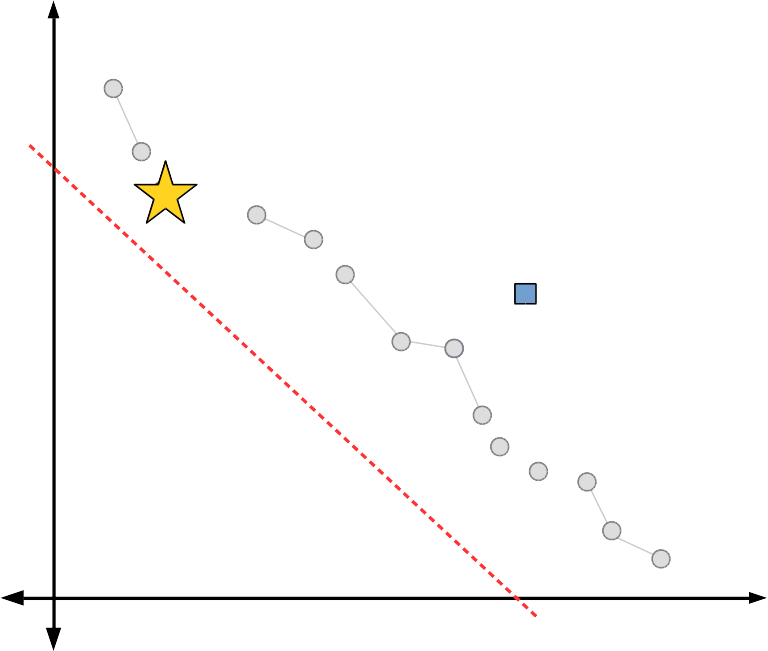}
\caption{Weighted sum method}
\label{preprocessing_weighted}
\end{subfigure}
\caption{Bound information when a single objective MILP terminates early. \revision{(i) $Y_N$, which we assume to be unknown, is shown in grey, (ii) the optimal solution, which we assume is not known at termination of the MILP solve, is depicted as a yellow star, (iii) the best known solution at termination is shown as a blue square, and (iv) the level curve associated with the best known dual bound at termination is shown as a dotted red line. For \ref{preprocessing_epsilon}, we assume that $\epsilon$ is defined so that the feasible region is restricted to the light blue box.} }
\label{preprocessing_fig}
\end{figure}

As can be seen, there is a certain level of trade-off present between the $\epsilon$-constraint method and the weighted sum method. The pros and cons of each technique are illustrated in Figures \ref{preprocessing_epsilon} and \ref{preprocessing_weighted}.

We now present Algorithms \ref{alg5} and \ref{alg6} in which we describe our proposed $\epsilon$-constraint and weighted sum based preprocessing procedures. On line 3 of Algorithm \ref{alg5} we solve the MILP associated with
\revision{\[
f_\lambda := (1 - \lambda) f_{1} + \lambda f_{2}.
\] 
The parameter} $\lambda$ is computed so that the level curves of $f_\lambda$ have the same slope as the line segment joining $y_I^1$ and $y_I^2$. On line 5 we then use the solution of this MILP to compute horizontal and vertical step sizes, $h_1$ and $h_2$, \revision{using the default value of $M = 60$.} These step sizes are then used to sequentially increase the values of $\epsilon_1$ and $\epsilon_2$ which are used on line 7 to construct new MILPs, using the $\epsilon$-constraint problem, which may yield new, undiscovered Pareto solutions. On lines 8 and 9 we modify the step sizes $h_1$ and $h_2$. If the MILP solved on line 7 yields a new, previously undiscovered Pareto solution, we decrease the step size. Otherwise, we increase it. This allows us the continue searching for additional new solutions in locations of $\OS$ which are near previously discovered solutions, and to cease searching in areas in which new solutions are not being generated. Note that the amount in which the step sizes are increased or decreased depends on the value of the parameter $\rho$. Also note that each time we solve a MILP, we utilize its solution to update $\N_s$. 

\begin{algorithm}[h!] 
  \caption{Preprocessing based on the $\epsilon$-constraint method.\\
  \underline{Input}: $y^1_I$, $y^2_I$ and a nonnegative value for parameter $\rho$.\\
  \underline{Output}: An initialized set of Pareto solutions $\protect\N_0 \subseteq Y_N$.}
  \label{alg5}
  \begin{algorithmic}[1]
    \Function{PreprocessingMethod1}{$y^1_I,y^2_I,\rho$} 
    	\State Let $\N_0 = \emptyset$.
    	\State Solve the MILP $\min\{f_\lambda(x): x\in {X}_I\}$ to obtain ${y}^\lambda_I \in {Y}_I$.
    	\State Add a cutting plane to $X$ lying on the level curve of $f_\lambda$ associated with the best dual solution.
    	\State Set $h_1 = \frac{(y^2_I)_1 - (y^\lambda_I)_1}{M}$, $\epsilon_1 = (y^\lambda_I)_1 + h_1$, $h_2 = \frac{(y^1_I)_2 - (y^\lambda_I)_2}{M}$ and $\epsilon_2 = (y^\lambda_I)_2 + h_2$, \revision{some $M \gg 1$}
    	\For{$k \in \{1,2\}$}{}
    		\While{$\epsilon_k > (y^k_I)_k$}{
			
			Solve the MILP $P_k(\epsilon_k):=\min\{f_{\{1,2\}\setminus\{k\}}(x): x\in {X}_I, f_k(x) \leq \epsilon_k\}$ to obtain $y^* \in Y_N$.}
    			\If{$\N_0 \not\dom y^*$}{ set $h_k = \frac{h_k}{1+\rho}$.}
    			\Else{ set $h_k = \max(5-\rho,1) h_k$.} 
    			\EndIf
    			\For{each $x \in X_I$ found while solving $P_k(\epsilon_k)$}{ 
			
			Let $N = $ \textsc{GenerateDualBd}$(s(x))$ and set $\N_0 = \nd{\N_0 \cup N}$.}
    			\EndFor
    			\State Set $\epsilon_k = \epsilon_k + h_k$.
    		\EndWhile
    	\EndFor
    	\State Return $\N_0$.
    \EndFunction
  \end{algorithmic}
\end{algorithm}

\begin{algorithm}[h!] 
  \caption{Preprocessing based on the weighted-sum method.\\
  \underline{Input}: A nonnegative value for parameter $\rho$.\\
  \underline{Output}: An initialized set of Pareto solutions $\protect\N_0 \subseteq Y_N$.}
  \label{alg6}
  \begin{algorithmic}[1]
    \Function{PreprocessingMethod2}{$\rho$} 
    	\State Let $\N_0 = \emptyset$.
    	\State Set $\Lambda = \{\lambda\}$, $\Lambda' = \{0,1\}$ and $t=0$.
    	\While{$t \leq \rho$}{}
    		\State Set $\tau = 0$ and $\sigma = |\Lambda|$.
    		\For{$\lambda' \in \Lambda$}{ remove $\lambda'$ from $\Lambda$ and add it to $\Lambda'$.\hfill(Assume $\Lambda'$ is always sorted in increasing order.)}
    			\State Solve the MILP $P(\lambda'):=\min\{f_{\lambda'}(x): x\in {X}_I\}$ to obtain ${y}^{\lambda'} \in {Y}_I$.
    			\State Add a cutting plane to $X$ lying on the level curve of $f_{\lambda'}$ associated with the best dual solution.
    			\If{$\N_0 \not\dom {y}^{\lambda'}$}{ set $\tau = \tau + 1$.}
    			\EndIf
    			\For{each $x \in X_I$ found while solving $P(\lambda')$}{ let $N = $ \textsc{GenerateDualBd}$(s(x))$ and set $\N_0 = \nd{\N_0 \cup N}$.}
    			\EndFor
    		\EndFor
    		\For{each adjacent pair $(\lambda_1,\lambda_2)\in\Lambda'$}{ add $\frac{\lambda_1+\lambda_2}{2}$ to $\Lambda$.}
    		\EndFor
    		\If{$\tau < \frac{\sigma}{5}$}{ set $t = t+1$.}
    		\EndIf
    	\EndWhile
    	\State Return $\N_0$.
    \EndFunction
  \end{algorithmic}
\end{algorithm} 

In Algorithm \ref{alg6} we compute several sets of weights which we utilize in the weighted-sum approach to generate Pareto solutions. We initialize the set of weights $\Lambda$ on line 3 with the weight $\lambda$ for which the level curves of $f_\lambda$ have the same slope as the line segment joining $y_I^1$ and $y_I^2$. We use $\sigma$ to represent the number of weights for which MILPs will be solved in a given iteration. We deem an iteration successful if at least a fifth of the solved MILPs reveal previously undiscovered Pareto solutions. We use $\tau$ to count the number of unsuccessful iterations. On line 11 we increase the number of weights that will be used in the next iteration by computing the next set of weights so that it contains the midpoint of each pair of adjacent weights in the set $\Lambda'$, which is the set of previously used weights together with 0 and 1. The process then terminates when the number of unsuccessful iterations exceeds the value of the parameter $\rho$. As we did with Algorithm \ref{alg5}, we also utilize the solution of each MILP we solve in this procedure to update $\N_s$.

\subsection{Probing} \label{sec:probe}
After Preprocessing, a probing technique can be used to strengthen the bounds on each integer variable, as stated below. 

\begin{proposition}[Probing on $x_{i}$]\label{probe}
Let $x_{i}$ be an integer variable. Fix $x_i = l_i$, relax integrality on other integer variables and solve the BOLP relaxation to obtain its Pareto set $\mathcal{L}_{l_i}$. If $\U_0 \dom \mathcal{L}_{l_i}$ then $X_{E}\subseteq\{x\colon x_{i}\ge l_{i}+1 \}$. 
\end{proposition}
\begin{proof}
Recognize that $\mathcal{L}_{l_i}$ dominates every $y \in Y_I$ where $y = \objf(x)$ with $x_i = \l_i$. The desired result follows from $\U_0 \dom \mathcal{L}_{l_i}$.
\end{proof}

This probing procedure can be repeated multiple times for a given integer $x_i$ and then iterated over each additional integer variable $x_j$. Furthermore, a similar procedure to that of Proposition \ref{probe} exists for tightening the upper bound. We point out that there are likely many more tasks that could be performed during Presolve and/or Preprocessing that could further impact the performance of BB. However, our goal here is not to develop extensive procedures for these tasks, but to put together an initial implementation that highlights some of what can be done.\\

\section{Node processing} \label{sec:np}
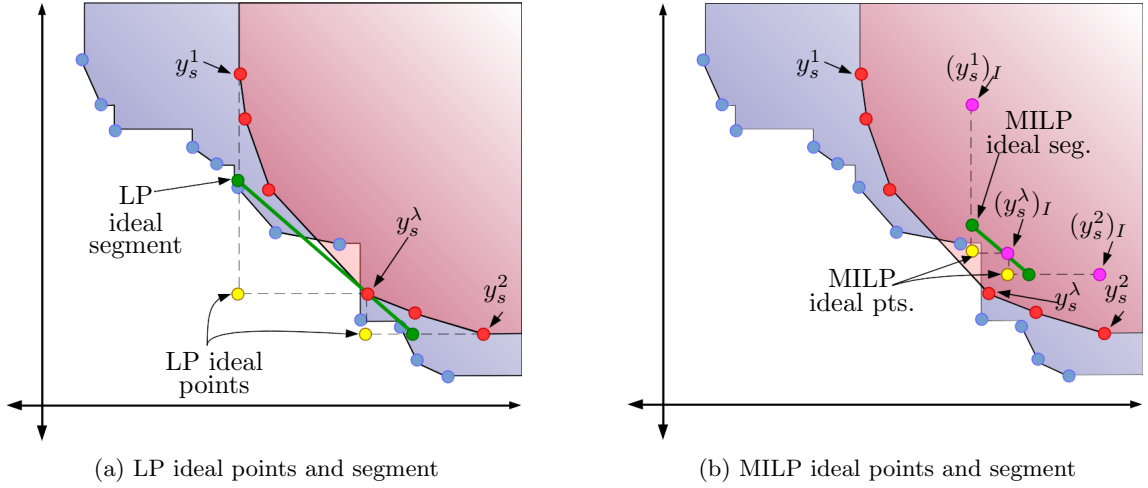
\begin{figure}
\begin{subfigure}[h]{.5\textwidth}
\centering
\begin{tikzpicture}
\node[anchor=south west,inner sep=0] (image) at (0,0) {\includegraphics[width=.83\textwidth]{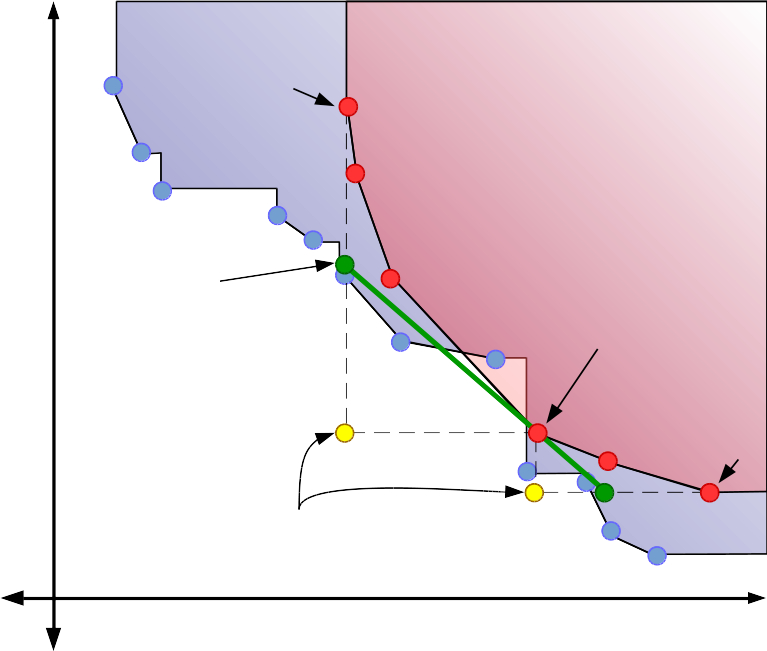}};
\node[align=center,yshift=2.1cm,xshift=-0.95cm] at (image.center) {
$y^1_s$ };
\node[align=center,yshift=0.0cm,xshift=2cm] at (image.center) {
$y^\lambda_s$ };
\node[align=center,yshift=-0.9cm,xshift=3.15cm] at (image.center) {
 $y^2_s$ };
\node[align=center,yshift=-1.85cm,xshift=-.6cm] at (image.center) {
 LP ideal };
 \node[align=center,yshift=-2.2cm,xshift=-.6cm] at (image.center) {
 points };
\node[align=center,yshift=.35cm,xshift=-1.7cm] at (image.center) {
LP} ;
\node[align=center,yshift=-.00cm,xshift=-1.7cm] at (image.center) {
ideal} ;
\node[align=center,yshift=-.35cm,xshift=-1.7cm] at (image.center) {
segment};
\end{tikzpicture}
\caption{LP ideal points and segment}
\label{ideal_pts_sgmnt}
\end{subfigure}%
\begin{subfigure}[h]{.5\textwidth}
\centering
\centering
\begin{tikzpicture}
\node[anchor=south west,inner sep=0] (image) at (0,0) {\includegraphics[width=.83\textwidth]{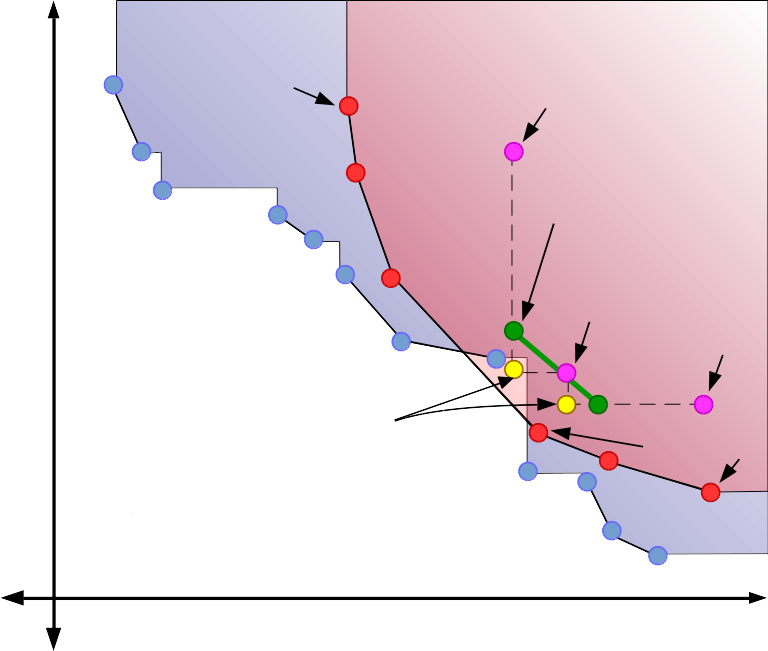}};
\node[align=center,yshift=2.1cm,xshift=-0.95cm] at (image.center) {
$y^1_s$ };
\node[align=center,yshift=-1.0cm,xshift=2.5cm] at (image.center) {
$y^\lambda_s$ };
\node[align=center,yshift=-0.9cm,xshift=3.15cm] at (image.center) {
 $y^2_s$ };
\node[align=center,yshift=-.75cm,xshift=-.3cm] at (image.center) {
 MILP} ;
\node[align=center,yshift=-1.1cm,xshift=-.3cm] at (image.center) {
ideal pts.} ;
\node[align=center,yshift=2.05cm,xshift=1.25cm] at (image.center) {
$(y^1_s)_I$ };
\node[align=center,yshift=0.25cm,xshift=1.9cm] at (image.center) {
$(y^\lambda_s)_I$ };
\node[align=center,yshift=-0.05cm,xshift=2.9cm] at (image.center) {
$(y^2_s)_I$ };
\node[align=center,yshift=1.35cm,xshift=2.05cm] at (image.center) {
MILP} ;
\node[align=center,yshift=1cm,xshift=2.05cm] at (image.center) {
ideal seg.} ;
\end{tikzpicture}
\caption{MILP ideal points and segment}
\label{ideal_pts_sgmnt2}
\end{subfigure}
\caption{Fathoming in biobjective BB}
\label{fath_rules_fig}
\end{figure}

Processing a node consists of three basic steps: (i) Generate a valid dual bound; (ii) Check a fathoming rule to determine whether or not $s$ can be eliminated from the search tree; (iii) Optionally, if $s$ is not fathomed in (ii), generate a tighter dual bound and repeat (ii). Figure \ref{fath_rules_fig} provides a visual example of how one might carry out these three steps. Most of the fathoming rules for biobjective BB are designed to check whether or not $\U_s$ dominates $(Y_s)_I$ by exploiting the transitivity of dominance. First, a set $\T$ is generated such that $\T \dom (Y_s)_I$. Then if $\U_s \dom \T$, $\U_s \dom (Y_s)_I$ and $s$ can be fathomed. Otherwise, a tighter bound on $(Y_s)_I$ is needed. The first bound we use is a set of two ideal points which we obtain by solving three single objective LPs; one for each $f_k$ and an one with a weighted sum objective $f_\lambda$ in which the weights, denoted $\lambda^s$, are given by the normal vector of the line segment $H_s$ passing through $y^1_s$ and $y^2_s$. We begin with these points because it is straightforward to determine whether or not $\U_s$ dominates a singleton. In Figure \ref{ideal_pts_sgmnt} these points are labelled ``LP ideal points.'' Notice that they are not dominated. Consider the intersection of $(Y_s)^{ideal} + \nonneg$ and the line with normal vector $\lambda^s$ passing through $y^\lambda_s$. Recognize that this intersection, which we denote $H_s^\lambda$, is also a valid dual bound. In Figure \ref{ideal_pts_sgmnt} the resulting line segment is labelled ``LP ideal segment,'' but is not dominated. A tighter bound can next be found by explicitly generating $\L_s$. In Figure \ref{ideal_pts_sgmnt} this is the set indicated by the red points, which is again not dominated. After generating $\L_s$, one cannot hope to find a tighter bound on $(Y_s)_I$ resulting from LP solutions. Instead, one can solve single objective MILPs to generate elements of $(Y_s)_I$ and use these elements to form a valid dual bound. We first generate ideal points in the same way as before, but use single objective MILPs rather than LPs. In Figure \ref{ideal_pts_sgmnt2} these points are labelled ``MILP ideal points.'' Yet again they are not dominated. We can then consider the intersection of $((Y_s)_I)^{ideal} + \nonneg$ and the line with normal vector $\lambda^s$ passing through $(y^\lambda_s)_I$, which we denote $\tilde{H}_s^\lambda$. This intersection forms another valid dual bound. In Figure \ref{ideal_pts_sgmnt2} the resulting line segment is labelled ``MILP ideal segment'' and is dominated. Hence, $s$ can be fathomed in this example.

We now formally outline the fathoming rules employed in this work. Some additional notation will be useful. For $k \in \{1,2\}$, define 
\begin{equation}
\mathcal{P}^k_s := \revision{\{y^{i}_{s} \colon i \neq k\}  \cup \{y^\lambda_s\},}
\end{equation} 
and let 
\begin{equation}
\mathcal{P}_s := (\mathcal{P}^1_s)^{ideal} \cup (\mathcal{P}^2_s)^{ideal}.
\end{equation} 
Additionally, for any $\mathcal{I} \subset \{1,2,\lambda\}$, define 
\begin{equation}
D^\mathcal{I}_s:=\cup_{k=1}^{2} \left(\left(\mathcal{P}^k_s\setminus \cup_{i\in \mathcal{I}}\, {y}^i_s \right) \cup \cup_{i \in \mathcal{I}\setminus\{k\}}\, ({y}^i_s)_I\right)^{ideal}.
\end{equation}
$\mathcal{P}_s$ represents the sets of ideal points obtained from LP solutions, while  $D^\mathcal{I}_s$ represents a set of ideal points obtained from a mixture of LP and MILP solutions. Our five fathoming rules are given below. Rule 0 expresses the idea of fathoming due to optimality, while the remainder of the rules indicate situations in which $s$ can be fathomed due to bound dominance.

\begin{proposition}[Fathoming Rules]\label{fath0} 
Node $s$ can be fathomed if any of the following holds:
\begin{enumerate}
\item[0.] $\L_s \subset (Y_s)_I$,
\item[1a.] $\U_s \dom \mathcal{P}_s$,
\item[2a.] $\U_s \dom {H}_s^\lambda$,
\item[1b.] $\U_s \dom D^\mathcal{I}_s$ for some $\mathcal{I} \subset \{1,2,\lambda\}$,
\item[2b.] $\U_s \dom \tilde{H}_s^\lambda$,
\item[3.] $\L_s \subseteq \U_s$.
\end{enumerate}
\end{proposition}
\begin{proof}
Rule 0 is due to integer feasibility of $\L_s$. Rule 1a holds since by construction $\mathcal{P}_s \dom \L_s$, and so $\U_s \dom \L_s$. Rule 2a holds since by construction $\tilde{H}_s^\lambda \dom \L_s$, and so $\U_s \dom \L_s$. For Rule 1b, note that by construction, for any $\mathcal{I} \subset \{1,2,\lambda\}$, $D^\mathcal{I}_s \dom (y_s)_I$ for every $(y_{s})_{I} \in (Y_s)_I$ and thus $D^\mathcal{I}_s$ is a valid dual bound at node $s$. For Rule 2b, note that by construction ${H}_s^\lambda \dom (y_s)_I$ for every $(y_s)_I \in (Y_s)_I$ and thus ${H}_s^\lambda$ is a valid dual bound at node $s$. Rule 3 is obvious. 
\end{proof}

Before we outline the process we use for processing a node $s$, we briefly discuss another important task that ought to be carried out while processing node $s$: Updating $\N_s$. We do this in two ways: (i) add each integer-feasible line segment discovered while checking Fathoming Rule 0 to $\N_s$, and (ii) for each discovered $x^* \in X_I$, generate the nondominated subset of 
\begin{equation}
\mathcal{Y}(x^*) := \{y = \objf(x): x \in X, x_i = x^*_i \text{ for } \revision{i = 1,\dots,n}\}
\end{equation}
and add each defining line segment of this set to $\N_s$. Consider the latter of these strategies. Observe that the feasible set of $\mathcal{Y}(x^*)$ can be interpreted as a leaf node of the BB tree, which we denote $s(x^*)$. Hence, the $\mathcal{Y}(x^*) + \nonneg = \L_{s(x^*)}$. This leads to a need for generating the nondominated subset of $\L_{s}$, i.e. $\nd{\L_s}$. Typical techniques for generating $\nd{\L_{s}}$ include the multiobjective simplex method and the parametric simplex algorithm (PSA) \cite{ehrgott2005multicriteria}. However, the multiobjective simplex method is far more robust than is necessary for biobjective problems. Also, we found in practice that using the PSA often resulted in many basis changes yielding the same extreme point of $\L_s$ in $\OS$. Since much work is done during the PSA to determine the entering and exiting variables, we found that generating $\nd{\L_s}$ using the PSA required a significant amount of computational effort. We decided to use an alternative method for generating $\nd{\L_s}$ which relies on sensitivity analysis. We first solve the single objective LP using objective $f_2$ to obtain $y^2_s$. Next we create the LP 
\begin{equation}
\mathscr{P}_s(\alpha):=\min\{ f_1(x) + \alpha f_2(x): x \in X_s\}
\end{equation}
and then carry out the procedure outlined in Algorithm \ref{alg2}.

\begin{algorithm}[h] 
  \caption{Generate $\nd{\L_s}$\\
  \underline{Input}: Node $s$.\\
  \underline{Output}: A set $\mathcal{B}$ containing all defining line segments of $\nd{\L_s}$.}
  \label{alg2}
  \begin{algorithmic}[1]
    \Function{GenerateDualBd}{$s$} 
        \State Set $\mathcal{B} = \emptyset$.
        \State Solve the LP $\min\{f_2(x): x \in X_s\}$ to obtain $y^2_s$.
    	\State Solve $\mathscr{P}_s(0)$ to obtain solution $x^*$ and set $y = \objf(x^*)$.
    	\While{$y \neq y^2_s$}{}
    		\State Use sensitivity analysis to obtain an interval $[\alpha',\alpha'']$ such that $x^*$ is optimal to $\mathscr{P}_s(\alpha)$ for all $\alpha \in [\alpha',\alpha'']$.
    		\State Let $\alpha^*$ be the negative reciprocal of the slope of the line through $y$ and $y^2_s$.
    		\State Set $x^* = \argmin\{\mathscr{P}_s(\alpha'' + \epsilon)\}$ for sufficiently small $\epsilon \in (0,\alpha^*-\alpha'']$.
    		\If{$\objf(x^*) \neq y$}
    			\State {Add the line segment connecting $\objf(x^*)$ and $y$ to $\mathcal{B}$. Update $y$ to be $\objf(x^*)$.}
    		\EndIf
    	\EndWhile
    	\State Return $\mathcal{B}$.
    \EndFunction
  \end{algorithmic}
\end{algorithm} 

In lines 3 and 4 of Algorithm \ref{alg2} we compute the south-east and north-west most extreme points of $\nd{\L_s}$, respectively. The while loop beginning on line 5 is then used to sequentially compute adjacent extreme points of $\nd{\L_s}$ in a west to east pattern, until the south-east most extreme point is rediscovered. Each line segment joining a pair of adjacent extreme points of $\nd{\L_s}$ is stored and the set of all computed segments is returned at the end of the procedure. Note that the correctness of the algorithm relies on an appropriately small choice for $\epsilon$ on line 8 (in our implementation, $\epsilon$ is set to $10^{-5}$). As we have discussed, there are other methods which can be used here that do not rely on $\epsilon$, such as the PSA or the first phase of the two-phase method for solving biobjective combinatorial problems \cite{ehrgott2005multicriteria}. We have already discussed the difficulties we encountered with the PSA. The difficulty with the first phase of the two-phase method is that, although it generates the extreme supported Pareto solutions of a BOLP, it does not generate them in order from left to right. Thus, when using a simplex-style solution method for each single objective LP, each iteration can require a significant number of basis changes. Our method generates these extreme points in order from left to right, and as a result, warm-starting each iteration by reusing the basis information from the previous iteration reduces the overall number of required basis changes.

Recognize from Proposition \ref{fath0} that Fathoming Rules 0 and 3 each impose a condition on $\L_s$ and therefore require knowledge of $\nd{\L_s}$ in order to be employed. We note, however, that for each of these rules it is often unnecessary to generate $\nd{\L_s}$ entirely. In particular, the generation of $\nd{\L_s}$ should cease if: (i) one is checking Fathoming Rule 0 and a defining line segment of $\nd{\L_s}$ is generated that is not integer feasible, or (ii) one is checking Fathoming Rule 3 and a defining line segment of $\nd{\L_s}$ is generated that is not contained in $\U_s$. Hence, the procedures in Algorithm \ref{alg2} can be modified in order to develop strategies for checking Fathoming Rules 0 and 3. These strategies are outlined in Algorithms \ref{alg4} and \ref{alg3}, respectively. 

\begin{algorithm}
  \caption{Fathoming Rule 0\\
  \underline{Input}: Node $s$ and solutions $y^1_s$ and $y^2_s$.\\
  \underline{Output}: 1 if node $s$ should be fathomed, 0 otherwise.}
  \label{alg4}
  \begin{algorithmic}[1]
    \Function{FR\_0}{$s, y^1_s, y^2_s$} 
    	\State $y^1_s$ is the solution to $\mathscr{P}_s(0)$. Let $x^*$ represent the preimage of $y^1_s$. Set $y = y^1_s$. 
    	\If{$y = y^2_s$}{ return 1}
    	\Else
    		
    	\While{$y \neq y^2_s$}{}
    		\State Use sensitivity analysis to obtain an interval $[\alpha',\alpha'']$ such that $x^*$ is optimal to $\mathscr{P}_s(\alpha)$ for all $\alpha \in [\alpha',\alpha'']$.
    		\State Let $\alpha^*$ be the negative reciprocal of the slope of the line through $y$ and $y^2_s$.
    		\State Set $x^* = \argmin\{\mathscr{P}_s(\alpha'' + \epsilon)\}$ for sufficiently small $\epsilon \in (0,\alpha^*-\alpha'']$.
    		\If{$\objf(x^*) \neq y$}
    			\State {Let $\mathbb{S}$ represent the line segment connecting $\objf(x^*)$ and $y$.}
    			\If{$\mathbb{S} \not\subset (Y_s)_I$}{ return 0}
    			\Else{ Update $y$ to be $\objf(x^*)$.}
    			\EndIf
    		\EndIf
    	\EndWhile
    	\State{return 1}
    	\EndIf
    \EndFunction
  \end{algorithmic}
\end{algorithm} 

Algorithm \ref{alg4} follows almost the same procedure as Algorithm \ref{alg2}, except it terminates prematurely on line 10 if a line segment is computed that is not integer feasible.
\begin{algorithm}
  \caption{Fathoming Rule 3\\
  \underline{Input}: Node $s$ and solutions $y^1_s$ and $y^2_s$.\\
  \underline{Output}: 1 if node $s$ should be fathomed, 0 otherwise.}
  \label{alg3}
  \begin{algorithmic}[1]
    \Function{FR\_3}{$s, y^1_s, y^2_s$} 
    	\State $y^1_s$ is the solution to $\mathscr{P}_s(0)$. Let $x^*$ represent the preimage of $y^1_s$. Set $y = y^1_s$. 
    	\If{$y = y^2_s$}
    		\If{$\U_s \dom y$}{ return 1}
    		\Else{ return 0}
    		\EndIf
    	\Else
    		
    	\While{$y \neq y^2_s$}{}
    		\State Use sensitivity analysis to obtain an interval $[\alpha',\alpha'']$ such that $x^*$ is optimal to $\mathscr{P}_s(\alpha)$ for all $\alpha \in [\alpha',\alpha'']$.
    		\State Let $\alpha^*$ be the negative reciprocal of the slope of the line through $y$ and $y^2_s$.
    		\State Set $x^* = \argmin\{\mathscr{P}_s(\alpha'' + \epsilon)\}$ for sufficiently small $\epsilon \in (0,\alpha^*-\alpha'']$.
    		\If{$\objf(x^*) \neq y$}
    			\State {Let $\mathbb{S}$ represent the line segment connecting $\objf(x^*)$ and $y$.}
    			\If{$\U_s \not\dom \mathbb{S}$}{ return 0}
    			\Else{ Update $y$ to be $\objf(x^*)$.}
    			\EndIf
    		\EndIf
    	\EndWhile
    	\State{return 1}
    	\EndIf
    \EndFunction
  \end{algorithmic}
\end{algorithm} 
Algorithm \ref{alg3} also follows almost the same procedure as Algorithm \ref{alg2}. However, this procedure terminates prematurely on line 5 or 12 if a point or line segment is computed that is not dominated by $\U_s$. We have now built the tools necessary to present our proposed procedure for processing a node $s$. We do so in Algorithm~\ref{alg1}.

\begin{algorithm}
  \caption{Process node $s$}
  \label{alg1}
  \begin{algorithmic}[1]
    \Function{ProcessNode}{$s$} 
    	\State Compute valid cutting planes for $(X_s)_I$ and add them to the description of $X_s$.
    	\For{$k \in \{1,2\}$}{ Solve $\min\{f_k(x): x\in {X}_s\}$ to find optimal solution $\bar{x}^k$ and generate ${y}^k_s \in {Y}^k_s$.}
    		\If{$y^k_s \in (Y_s)_I$}{ let $N = $ \textsc{GenerateDualBd}$(s(\bar{x}^k))$ and set $\N_s = \nd{\N_s \cup N}$.}
    		\EndIf
    	\EndFor
    	\If{$y^1_s,y^2_s \in (Y_s)_I$}
    		\If{\textsc{FR\_0}$(s,y^1_s,y^2_s)$ = 1}{ Fathom $s$, STOP!\hfill(Fathoming Rule 0)}
    		\EndIf
    	\EndIf
    	\State{Calculate $H_s$ and $\lambda^s$ using ${y}^1_s$ and $y^2_s$. Solve $\min\{f_\lambda(x): x\in {X}_s\}$ to find optimal solution $\bar{x}^\lambda$ and generate ${y}^\lambda_s \in {Y}^\lambda_s$.}
    	\If{$y^\lambda_s \in (Y_s)_I$}{ let $N = $ \textsc{GenerateDualBd}$(s(\bar{x}^\lambda))$ and set $\N_s = \nd{\N_s \cup N}$.}
    	\EndIf
    	\If{$\U_s \dom {y}^1_s$, $\U_s \dom {y}^2_s$ and $\U_s \dom {y}^\lambda_s$}
    		\If{$\U_s \dom \mathcal{P}_s$}{ Fathom $s$, STOP! \hfill(Fathoming Rule 1a)}
    		\Else
    			\State{Calculate $\tilde{H}_s^\lambda$.}
    			\If{$\U_s \dom \tilde{H}_s^\lambda$}{ Fathom $s$, STOP! \hfill(Fathoming Rule 2a)}
    			\Else
					\If{\textsc{FR\_3}$(s,y^1_s,y^2_s)$ = 1}{ Fathom $s$, STOP!\hfill(Fathoming Rule 3)}
    				\EndIf
    			\EndIf
    		\EndIf
    	\Else
    		\State Define the set $\mathcal{I} = \emptyset$.
    		\For{$k \in \{1,2\}$}
    			\If{$\U_s \not\dom (\mathcal{P}^k_s)^{ideal}$}{ add $(\{1,2\}\setminus\{k\})\cup\{\lambda\}$ to $\mathcal{I}$}
    			\EndIf
    		\EndFor
    		\For{each $k \in \mathcal{I}$}{ solve the MILP $\min\{f_k(x): x\in ({X}_s)_I\}$ to find optimal solution $\hat{x}^k$ and obtain $({y}^k_s)_I \in ({Y}^k_s)_I$.}
    			\State Add a local cut to $X_s$ lying on the level curve of $f_k$ associated with the best dual solution.
    			\State Let $N = $ \textsc{GenerateDualBd}$(s(\hat{x}^k))$ and set $\N_s = \nd{\N_s \cup N}$.
    		\EndFor
			\If{$\U_s \dom D^\mathcal{I}_s$}{ Fathom $s$, STOP! \hfill(Fathoming Rule 1b)}
			\ElsIf{$\lambda \in \mathcal{I}$}
				\State{Calculate ${H}_s^\lambda$.}
    			\If{$\U_s \dom {H}_s^\lambda$}{ Fathom $s$, STOP! \hfill(Fathoming Rule 2b)}
    			\EndIf
			\EndIf    		
    	\EndIf
    \EndFunction
  \end{algorithmic}
\end{algorithm} 

Line 2 of Algorithm~\ref{alg1} is an optional procedure in which we can generate locally valid cutting planes to strengthen the representation of $X_s$ if so desired. We then compute $y_s^1$ and $y_s^2$ on line 3. We then check to see if either of these solutions are integer feasible, and if they are, we generate the dual bound associated with the integer solution in order to update $\N_s$. Furthermore, if both solutions are integer feasible, we check Fathoming Rule 0 on line 6. On line 7 we compute the value $\lambda_s$, the value of the weights on the objectives so that the level curves of $f_\lambda$ have the same slope as the line segment joining $y_s^1$ and $y_s^2$. We then solve the LP associated with $f_\lambda$. If the solution is integer feasible, we again update $\N_s$ as before. On line 9 we check whether or not $y_s^1,y_s^2$ and $y_s^\lambda$ are dominated by $\U_s$. If they are, we proceed to check Fathoming Rules 1a, 2a, and 3. Otherwise, we solve the MILP associated with $f_\lambda$ and $f_k$ for each $k\in\{1,2\}$ such that the ideal point $(\mathcal{P}_s^k)^{ideal}$ is not dominated by $\U_s$. On lines 21 and 22 we utilize the solutions of each MILP to (optionally) add local cuts to $X_s$ and update $\N_s$. Finally, we check Fathoming Rules 1b and 2b.

Two additional tasks are performed while processing each node. 

\subsection{Objective space fathoming}
After processing each node, we perform an additional type of fathoming which we refer to as \emph{objective-space fathoming}. After updating $\N_s$, we impose bounds on $f_1$ and $f_2$ which ``cut off'' portions of $\OS$ in which we have discovered that $\U_s \dom (Y_s)_I$. In certain cases the remaining subset of $\OS$ 
consists of disjoint regions. When this is the case, we implement objective-space fathoming by branching on $f_1$ and $f_2$ bounds which generate the desired disjunctions in $\OS$. In these cases, objective-space fathoming resembles the ``Pareto branching'' of \citet{stidsen2014branch}~ and ``objective branching'' of \citet{parragh2015boip}. 

\subsection{Bound tightening}\label{sec:tight}
\begin{figure}
\begin{subfigure}[h]{.32\textwidth}
\centering
\includegraphics[width=0.97\textwidth]{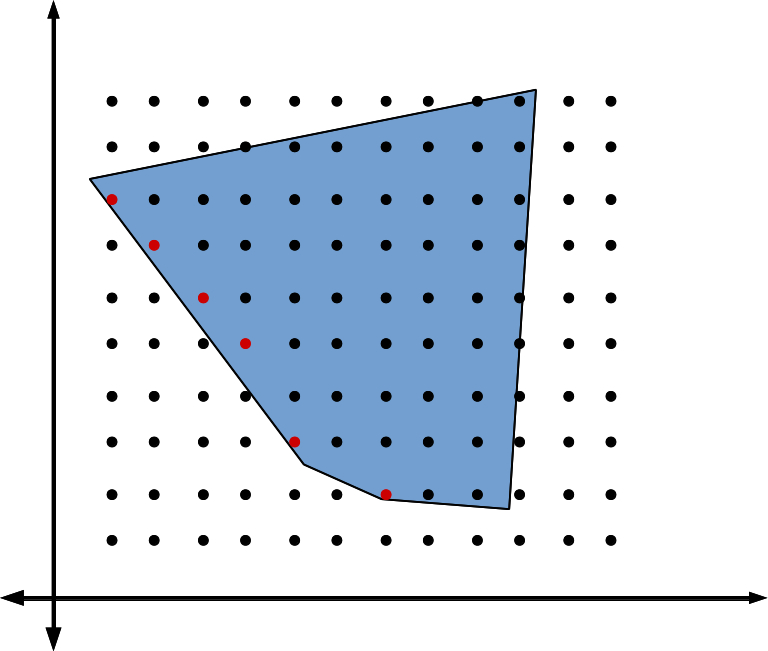}
\caption{Example instance of BOMILP}
\label{local_cuts_1}
\end{subfigure} 
\begin{subfigure}[h]{.33\textwidth}
\centering
\includegraphics[width=.97\textwidth]{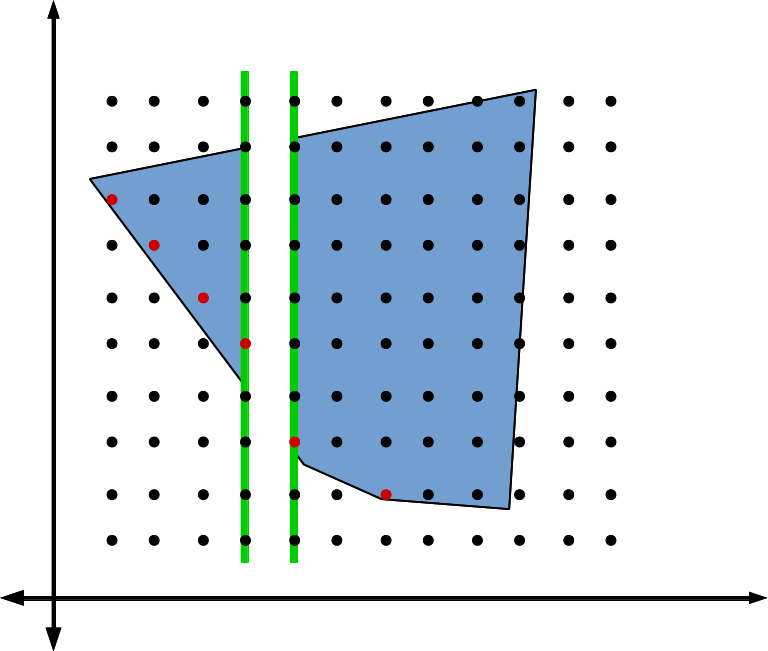}
\caption{After branching}
\label{local_cuts_2}
\end{subfigure}
\begin{subfigure}[h]{.32\textwidth}
\centering
\includegraphics[width=.97\textwidth]{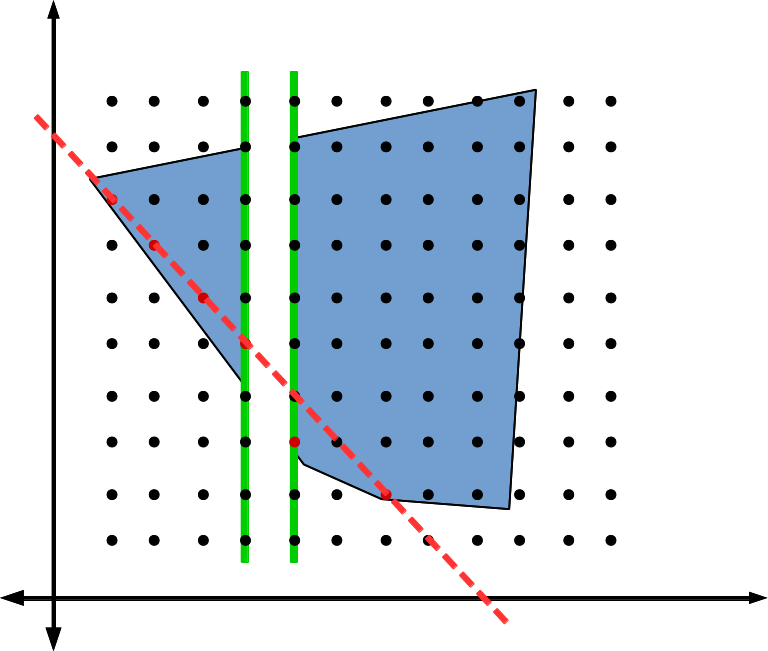}
\caption{Locally valid cut}
\label{local_cuts_3}
\end{subfigure}
\caption{An example showing the usefulness of locally valid cuts for BOMILP}
\label{local_cuts}
\end{figure}

In order to increase the likelihood of fathoming, we utilize a few different strategies for tightening the bound $\L_s$. The first strategy we use is the generation of locally valid cutting planes. We do this in two ways: (i) we generate discjuntive cuts based on disjunctions observed in $\OS$ when performing $\OS$ fathoming, and (ii) we convert the BOLP relaxation associated with $s$ to the BOMILP $\min\{f_\lambda(x): x\in ({X}_s)_I\}$, allow the MILP solver to process its root node, and add all cuts generated by this solver as local cuts to $s$ as local cuts. It is widely accepted that for single objective MILPs, locally valid cutting planes are not particularly helpful for improving the performance of BB. However, locally valid cutting planes can have a significantly greater impact on BOMILPs. To see this, observe Figure \ref{local_cuts}. Assume that Figure \ref{local_cuts_1} displays an instance of BOMILP for which the $(f_1,f_2)$-space and the $X$-space are one and the same, i.e., this instance contains only two variables $y_1$ and $y_2$, both integer, and $f_1 = y_1$ and $f_2 = y_2$. The constraints of this instance yield the blue polytope, and the integer lattice is indicated by the black dots. The red dots represent the Pareto-optimal solutions. Suppose that branching is performed as shown in Figure \ref{local_cuts_2}. Notice that all Pareto optimal solutions in the left branch can be revealed by a single locally valid cutting plane, as shown by the red dashed line in Figure \ref{local_cuts_3}. Also notice that this could never be accomplished through the use of globally valid cuts.

\subsection{Comparison with another \bb{}} \label{sec:compare}

We highlight some key differences regarding the node processing step between our \bb{} and that of \citet{belotti2012biobjective,belotti2015fathoming}, which is the only other \bb{} method for general BOMILP. There are also differences in the other components of \bb{}, 
but that is not of concern here.

The two methods differ in the way fathoming rules are implemented. Firstly, we utilize the data structure of \citet{treestructure} to store and dynamically update the set $\N_s$ throughout the BB process. In \cite{belotti2012biobjective,belotti2015fathoming}, fathoming rules are checked at a node $s$ of the BB tree by: (i) using $\N_s$ to generate $\U_s$ by adding a set of local nadir points to $\N_s$, (ii) selecting the  subset $\mathcal{R} := \U_{s}\cap ((Y_s)^{ideal} + \nonneg)$, and (iii) solving auxiliary LPs to determine whether $\mathcal{R}$ and $\L_{s}$ can be separated by a hyperplane. Node $s$ is then fathomed if $\mathcal{R} = \emptyset$ or if a separating hyperplane is found. Note that these procedures amount to comparing each element of the primal bound with the dual bound as a whole by solving at most one LP for each element of the primal bound.

In this paper, we utilize the opposite approach to fathoming. Rather than comparing each element of the primal bound with the dual bound as a whole, we compare each element of the dual bound with the primal bound as a whole. Additionally, instead of making these comparisons by solving LPs, we exploit the following guarantee of the data structure of \cite{treestructure}: a point or line segment inserted to the structure is added to the structure if and only if the point or segment is not dominated by the data already stored in the structure. Hence, we implement an extra function \textsc{IsDominated}($\cdot$) alongside this data structure which returns 1 if the input is dominated by $\N_s$ and 0 otherwise. We then implement our fathoming rules 1-3 by passing the appropriate sets ($\mathcal{P}_s, H_s^\lambda, D_s^\mathcal{I}, \tilde{H}_s^\lambda$ and $\L_s$) to \textsc{IsDominated}. If a 1 is returned for any of these sets, we fathom, otherwise we do not. It is difficult to comment on whether solving LPs or utilizing a function call to a data structure is more efficient for checking fathoming. However, we have found in practice that for a particular node $s$ of the BB tree, the primal bound $\U_s$ typically contains far more points and segments than the dual bound $\L_s$. Thus, comparing each element of the dual bound with the primal bound as a set seems to be a more efficient procedure than doing it the opposite way.

\section{Biobjective BB}\label{sec:multi_BB}
We now discuss the extension of the remaining major aspects of \bb{} to the biobjective setting. In particular, we discuss the specifics of how the different components of single objective BB ---  presolve/preprocessing, node processing, and branching, can each be extended to the biobjective setting. We then briefly discuss optional additions to our basic biobjective BB procedure. 

\subsection{Branching}\label{sec:branch}
In general, any rule for selecting a branching variable is permissible. However, it should be noted that for BOMILP several $y \in {Y}$, and consequently several $x\in {X}$, may be discovered while processing a node $s$. In fact, our implementation requires solving at least three LPs at each node. Since the variables may take on different values at each solution, it is possible that an integer variable takes a fractional value at some of these solutions and not at others. Because of this, we use a scoring scheme for branching in which each integer variable is given a score. Of the variables with the highest score, the one with the highest index is selected for branching. The score of $x_i$ is increased if: (i) $x_i$ is fractional at the LP solution associated with objective $f^k$, $k \in \{1,2,\lambda^s\}$, (ii) $x_i$ changes value at a pivoting step of Algorithm \ref{alg4}, or (iii) multiple single objective MILPs are solved to optimality at $s$ and $x_i$ takes different values for at least two of the MILP solutions.

After a branching decision has been made we utilize probing, as introduced in Proposition \ref{probe}, to strengthen bounds on each variable for both of the resulting subproblems. We do this for several reasons: (i) we may find during this process that our branching decision results in an infeasible subproblem, in which case we can discard the infeasible subproblem, enforce that the variable bounds associated with the feasible subproblem be satisfied at any child node of $s$, and choose a new branching variable; (ii) because much work in biobjective BB is dedicated to fathoming, we want to generate the strongest dual bound possible, which probing helps us to do; (iii) since processing a node in biobjective BB is an expensive operation, we seek to limit the number of nodes explored and probing aids in this endeavor by reducing the number of possible future branching decisions. We found during testing that this probing scheme at each node was extremely powerful, both in reducing the number of nodes processed during BB as well as overall running time. See Table \ref{table_presolve} in Section \ref{sec:compute} for evidence of this.

\subsection{Exploiting gaps in $\OS$}\label{sec:osgap}

\begin{figure}
\begin{subfigure}[h]{.32\textwidth}
\centering
\includegraphics[width=\textwidth]{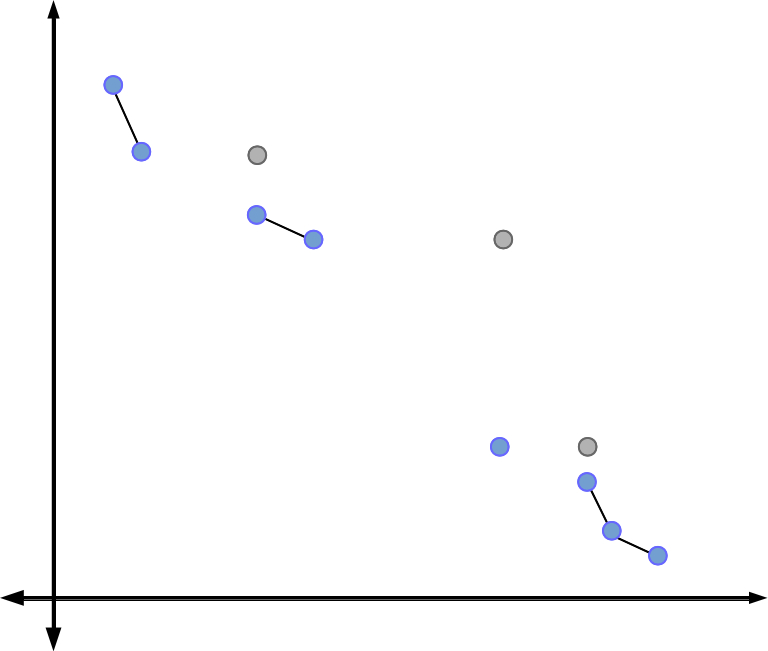}
\caption{Gaps}
\label{gaps1}
\end{subfigure} 
\begin{subfigure}[h]{.32\textwidth}
\centering
\includegraphics[width=\textwidth]{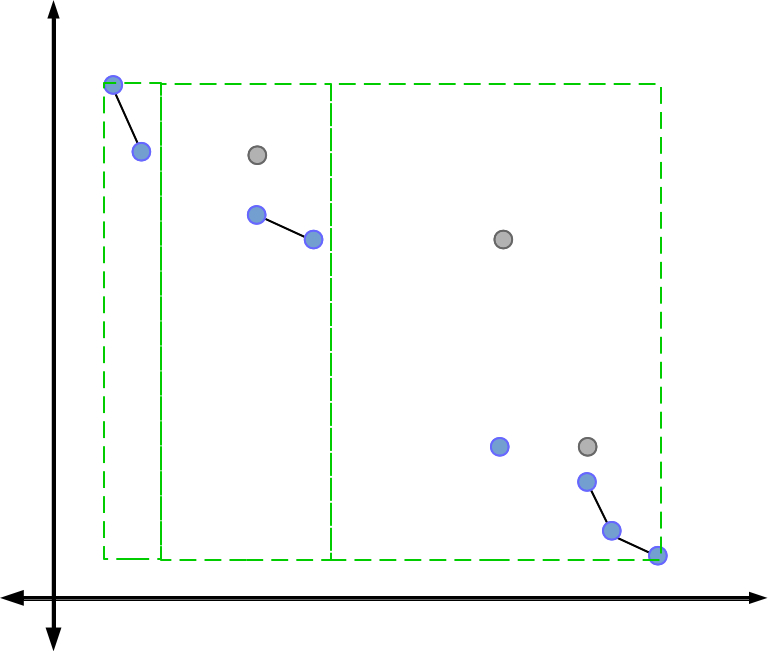}
\caption{Slitting $\OS$}
\label{gaps2}
\end{subfigure}
\begin{subfigure}[h]{.32\textwidth}
\centering
\includegraphics[width=\textwidth]{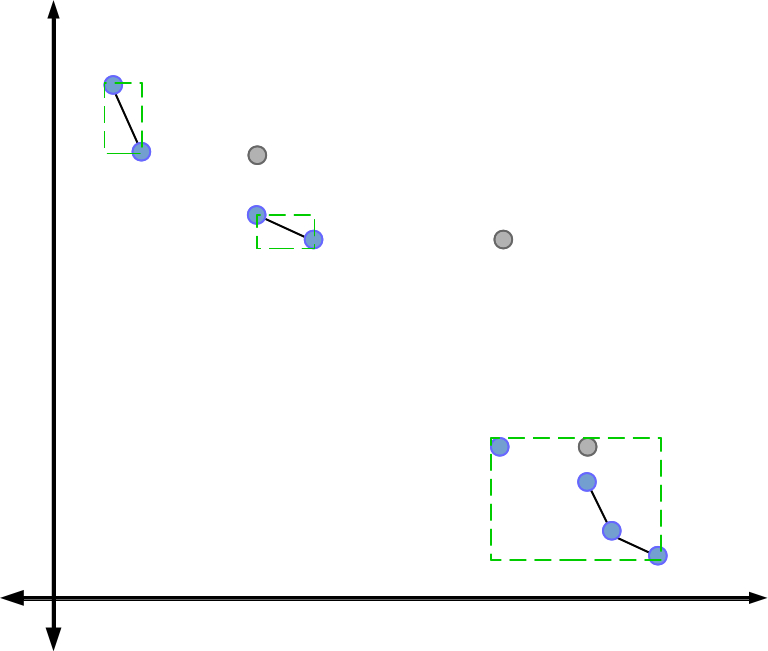}
\caption{Reducing the subregions}
\label{gaps3}
\end{subfigure}
\caption{Large gaps between solutions in $\OS$}
\label{gaps_fig}
\end{figure}

Due to the noncontinuous, nonconvex nature of the Pareto set of a BOMILP, there are occasionally large gaps between Pareto solutions in $\OS$. If this occurs, the likelihood that $\L_s \subseteq \U_s$ is significantly decreased for each node. Hence, this can result in an extreme amount of computational effort which yields no additional Pareto solutions. One way to combat this issue is to observe the solutions obtained during Preprocessing and record locations in $\OS$ where large gaps exist between discovered solutions. One can then split $\OS$ into a series of subregions based on the locations of these gaps and solve single objective MILPs (using objectives $f_1$ and $f_2$) within each subregion in order to remove locations containing no Pareto solutions. Afterwards BB can be run in each subregion rather than over the entire $\OS$. To aid in understanding this idea, observe Figure \ref{gaps_fig}. Here Pareto solutions are shown in blue and subregions in $\OS$ are indicated by green dashed lines.

\subsection{Measuring Performance} \label{sec:gap}

In single objective BB, one can terminate the procedure at any time and obtain a measure of the quality of the best known solution in terms of the gap between this solution and the best known dual bound. 
We propose a similar scheme for biobjective BB. 
Let $\mathscr{O}_{s^*}$ represent the set of open nodes after a node $s^*$ has been processed. After processing $s^*$, the global dual bound is $\mathcal{DB}_{s^*}=\nd{\cup_{s\in\mathscr{O}_{s^*}}\L_s}$. Therefore, if BB is terminated after $s^*$ is processed, the performance of BB can be quantified by measuring the distance between $\mathcal{DB}_{s^*}$ and $\U_{s^*}$. One natural metric to use for measuring this distance 
is the Hausdorff metric: \[d_H(\mathcal{DB}_{s^*},\U_{s^*}) := \max \left\{\sup_{i\in \mathcal{DB}_{s^*}}\inf_{j\in \U_{s^*}} d(i,j),\, \sup_{j\in \U_{s^*}}\inf_{i\in \mathcal{DB}_{s^*}} d(i,j) \right\}.\] 
Unfortunately the nonconvex nature of $\U_s$ makes the Hausdorff metric difficult to use 
since it cannot be computed using a linear program. In our implementation $\U_{s^*}$ is stored as the individual line segments and singletons comprising $\N_{s^*}$ using the data structure of \citep{treestructure}. $\mathcal{DB}_{s^*}$ is computed by generating the points and line segments comprising its nondominated subset, which are also stored using the same data structure. Thus, rather than explicitly computing $d_H(\mathcal{DB}_{s^*},\U_{s^*})$, we instead compute 
\begin{equation*}
\mathcal{G}_{s^*}:=\max\{d_H(\mathcal{DB}_{s^*},\mathcal{S}+\nonneg): \mathcal{S} \in \N_{s^*}\}
\end{equation*}
via pairwise comparison of the points and line segments comprising $\mathcal{DB}_{s^*}$ and $\N_{s^*}$. Clearly, $\mathcal{G}_{s^*}$ is a upper bound on $d_H(\mathcal{DB}_{s^*},\U_{s^*})$. Recognize, though, that $\mathcal{G}_{s^*}$ is an absolute measurement and so it is difficult to use to compare the performance of BB on multiple instances of BOMILP. Thus, in practice we use a percentage calculated as 
\begin{equation*}
\bar{\mathcal{G}}_{s^*} := 100\times\frac{\left|\max\{y^2_1-y^1_1,y^1_2-y^2_2\} - \mathcal{G}_{s^*}\right|}{\max\{y^2_1-y^1_1,y^1_2-y^2_2\}}.
\end{equation*}
\revision{We refer to this number as the \% duality gap. The number $\bar{\mathcal{G}}_{0}$ means the duality gap at the root node.}

Another method for measuring the distance between $\mathcal{DB}_{s^*}$ and $\U_{s^*}$ is to compute a so called \emph{hypervolume gap}. Let $hv(\cdot)$ denote the area of subset of $\R^2$. Then the hypervolume gap between $\mathcal{DB}_{s^*}$ and $\U_{s^*}$, as proposed by \citet{zitzler2003performance}, is
\begin{equation*}
\mathcal{HV}_{s^*} := 100\times\frac{hv((\mathcal{DB}_{s^*} + \nonneg) \cap \OS) - hv(\U_{s^*} \cap \OS)}{hv((\mathcal{DB}_{s^*} + \nonneg) \cap \OS)},
\end{equation*} 
A similar measure is used to assess the quality of approximations to the Pareto sets of BOMILP instances in \citep{boland2015acriterion}. 

Recognize that the Hausdorff and hypervolume gap measurements play significantly different roles. The hypervolume gap provides a measure of the proximity of the dual bound to the primal bound throughout the entirety of $\OS$, while the Hausdorff gap provides a measure of the proximity of the dual and primal bounds in the location at which they are furthest apart. Hence, we can interpret the Hausdorff gap as a worst-case measurement and the hypervolume gap as a sort of average-case measurement. We note that in our initial tests we utilize both the Hausdorff and hypervolume measurements so that our results can be compared with other works, such as \citep{boland2015acriterion}, which use the hypervolume gap. However, since the Hausdorff gap provides a worst-case measure and is therefore more robust, we do not use the hypervolume gap measurement in our final set of experiments.

\subsection{\revision{Our BB algorithm}}

A pseudocode of our \bb{} procedure is given in Algorithm \ref{alg7}. 

\begin{algorithm}[h!] 
  \caption{\bb{} for BOMILP.\\
  \underline{Input}: An instance $\mathcal{I}$ of BOMILP.\\
  \underline{Output}: The Pareto set of instance $\mathcal{I}$.}
  \label{alg7}
  \begin{algorithmic}[1]
    \Function{BBsolve}{$\mathcal{I}$} 
    	\State Set $\mathscr{L} = \emptyset$.
    	\State Use primal presolve, biobjective duality fixing and exploitation of singleton and dominating columns to simplify $\mathcal{I}$.
    	\For{$k \in \{1,2\}$}{ solve the MILP $\min\{f_k(x): x\in {X}_I\}$ to obtain ${y}^k_I \in {Y}_I$.}
    	\EndFor
    	\State Select $\rho \geq 0$ and run either \textsc{PreprocessingMethod1}$(y^1_I,y^2_I,\rho)$ or \textsc{PreprocessingMethod2}$(y^1_I,y^2_I,\rho)$ to return $\N_0$.
    	\State Perform probing to further simplify $\mathcal{I}$.
    	\State Add the continuous relaxation of $\mathcal{I}$ to $\mathscr{L}$.
    	\While{$\mathscr{L} \neq \emptyset$}{ select $s$ from $\mathscr{L}$.}
    		\State Run \textsc{ProcessNode}$(s)$.
    		\If{$s$ is not fathomed}{ perform $\OS$ fathoming.}
    			\If{the nondominated portion of $\OS$ consists of disjoint regions}{ perform Pareto branching. Add the resulting subproblems to $\mathscr{L}$.}
    			\Else{ select the variable with highest score for branching.}
    				\State Perform probing to simplify each of the subproblems resulting from the current branching decision.
    				\If{probing reveals an infeasible subproblem}{ impose the restrictions of the feasible subproblem and select the variable with the next highest score for branching. Repeat Line 13.}
    				\Else{ branch on the selected variable. Add the resulting subproblems to $\mathscr{L}$.}
    				\EndIf
    			\EndIf 
    		\EndIf
    	\EndWhile
    	\State Return $\N_{s^*}$, where $s^*$ is the last node for which \textsc{ProcessNode} was called.
    \EndFunction
  \end{algorithmic}
\end{algorithm}

\section{Computational Analysis} \label{sec:compute}

We implemented Algorithm~\ref{alg7} for our \bb{} scheme using the C programming language and the ILOG CPLEX 12.6 optimization package. This implementation, along with the instances we generated for use in Section \ref{sec:MIPLIB} can be found at \url{https://github.com/nadelgr/BOMILP_BB}. \citet{boland2015acriterion} graciously shared their code with us and so we were able to run their triangle splitting (TS) algorithm on our machines to compare against our BB algorithm. Recall that TS is a search method in the objective space. In preliminary tests, we also compared with the BB method of \citep{belotti2012biobjective}. However, their implementation was incomplete and so the performance of our \bb{} was far superior to theirs. For this reason, we do not include the results of their \bb{}. All testing described in Sections \ref{sec:presolve_computation}--\ref{sec:approximation_computation} was conducted using a Dell PowerEdge R430 server running Fedora Core 27 and which had a Xeon E5-2640 CPU and 64 GB of RAM. \revision{For tests described in Section \ref{sec:MIPLIB} we utilized the Extreme Science and Engineering Discovery Environment (XSEDE) \citep{xsede} Bridges system at the Pittsburgh Supercomputing Center (PSC) through allocation DMS200019. Specifically, these tests were conducted using a HPE Apollo 2000 server running CentOS Linux 7 and which had a Intel Haswell CPU and 128 GB of RAM.}

For experiments described in \textsection\ref{sec:presolve_computation}--\textsection\ref{sec:approximation_computation} we utilized a test set consisting of the instances examined in \citet{belotti2012biobjective} and \citet{boland2015acriterion,boland2013criterion}. The former contained \revision{30 instances} with 60 variables and 60 constraints (\instance{Belotti60}) and \revision{30 instances} with 80 variables and 80 constraints (\instance{Belotti80}). The latter had \revision{5 instances} for each of the three types \instance{Boland80}, \instance{Boland160}, and \instance{Boland320} (we do not solve instances with less than 60 constraints or variables due to their relative ease), and \revision{4 instances} for each of the three types \instance{Boland16}, \instance{Boland25}, and \instance{Boland50}.\footnote{\revision{These are labelled this way to maintain consistency with the way other instances are labeled although the respective total number of variables and constraints is approximately 800, 1250 and 2500.}}

\begin{figure}[h!]
\begin{subfigure}[h]{.55\textwidth}
\centering
\includegraphics[width=.9\textwidth]{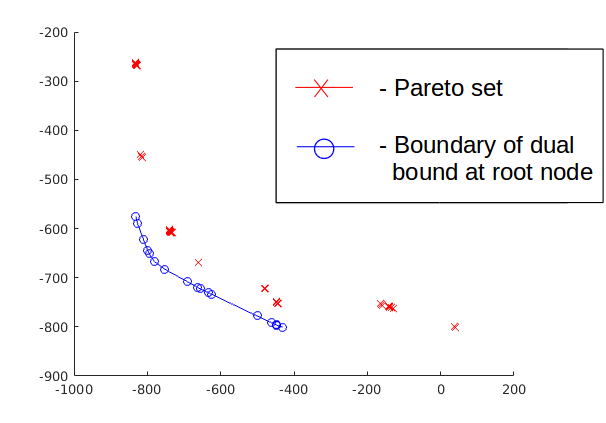}
\caption{Instance from the \instance{Belotti60} set.}
\label{belotti_instance}
\end{subfigure}%
\begin{subfigure}[h]{.45\textwidth}
\centering
\includegraphics[width=1\textwidth]{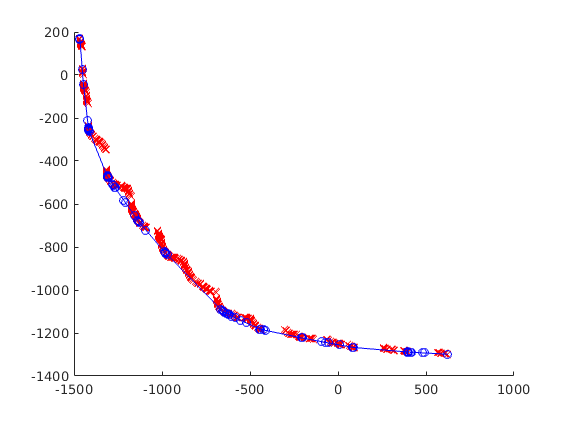}
\caption{Instance from the \instance{\instance{Boland16}} set.}
\label{boland_instance}
\end{subfigure}
\caption{Pareto set and boundary of $\L_0$ for the two instance families.}
\label{pareto_sets}
\end{figure}

\revision{
Figure \ref{pareto_sets} depicts the Pareto set and boundary of $\L_0$ for one instance from each of the two instance classes. Note the following structural differences displayed in the two figures
\begin{enumerate}
\item The relative gap between the Pareto set and boundary of $\L_0$ is greater in Figure \ref{belotti_instance} than in Figure \ref{boland_instance}.
\item The relative gap between connected subsets of the Pareto set is greater in Figure \ref{belotti_instance} than in Figure \ref{boland_instance}.
\item The overall number of solutions present in the Pareto set is greater in Figure \ref{boland_instance} than in Figure \ref{belotti_instance}.
\end{enumerate} 
We found that the above differences were typical for these instance families. This provides some insight into the differences in performance seen for these two instances families through the rest of this section. Note, in particular, the difference in duality gaps seen in Experiments 1, 3, and 5 as well as the difference in number of nodes processed when utilizing $\OS$ gap splitting in Experiment 4. 

Our final set of experiments are described in \textsection\ref{sec:MIPLIB}, where we opted to generate a more difficult test set. For this purpose, we created biobjective variants of instances from MIPLIB 2017 \citep{miplib2017pre}  that were feasible, mixed-integer, marked easy, and contained at most 1000 decision variables. For each such instance, we generated two secondary objective functions and discarded instances for which: (i) the Pareto set was a singleton, or (ii) the second objective was unbounded, or (iii) the MILP associated with either $f_1$ or $f_2$ took over 8 hours to solve.}

\revision{The computational tests with our BB had a maximum solution time of 8 hours. For each instance, we recorded the computation time in seconds, the number of nodes explored in our BB tree, and the \% duality gap computed after the root node was processed or at termination, as applicable. We report average values of these numbers for the \instance{Belotti*} instances, which we recall are 30 of each type, since the performance of these instances was mostly similar, and individual numbers for the other instances.}

We began our tests by turning off all nonessential features of our BB procedure, and then sequentially turning on various features to test their impact on the overall procedure. If a particular feature of our BB procedure was deemed effective in reducing the overall effort required to solve instances of BOMILP, this feature was left on for the remainder of the tests, otherwise it was turned back off. 

Our original implementation included a variety of features which did not prove useful in either reducing the overall BB time or the number of explored nodes. For the sake of space, in the sections that follow we focus only on features that proved useful. We briefly note some of these ideas here to motivate future research into them. Most of our fruitless features involved adding various cutting planes to the problem formulation. Note that we are not referring to CPLEX default cut generation -- this was left on and did prove useful. Instead, we are referring to: (i) attempts to add user-generated cuts from discovered disjunctions, and (ii) attempts to use CPLEX default cut generation \emph{at each node} and add the discovered cuts as local cuts. Other attempted features included checks for early termination of Fathoming Rule 3 and the generation of $\nd{\L_s}$. Each of these provided inconsistent results, reducing BB time for some problems but increasing it for others. Hence, both were abandoned in the end.

\subsection{Presolve Techniques}\label{sec:presolve_computation}

Table \ref{table_presolve} contains the results of our first computational experiment.  For this test we utilized \textsc{PreprocessingMethod2} with $\rho$ set to zero. \revision{The column ``Fixed'' in this table refers to the number of variables that were fixed by the presolve method.}

\begin{table}[h!]
\footnotesize
\centering
\caption{Experiment 1 -- Impact of presolve techniques}\label{table_presolve}
\begin{tabular}{l|rrr|rrrr|rrr|rrrr}
\toprule
\multirow{2}{*}{Instance} & \multicolumn{3}{c|}{All Off} & \multicolumn{4}{c|}{Duality Fixing} & \multicolumn{3}{c|}{Singleton Columns} & \multicolumn{4}{c}{Dominating Columns}\\
& Time & Nodes & $\bar{\mathcal{G}}_{0}$ & Time & Nodes & Fixed & $\bar{\mathcal{G}}_{0}$ & Time & Nodes & $\bar{\mathcal{G}}_{0}$ & Time & Nodes & Fixed & $\bar{\mathcal{G}}_{0}$  \\ 
\midrule
\instance{Belotti60} (30)	&	4	&	77	&	53	&	4	&	77	&	0	&	53	&	4	&	77	&	53	&	4	&	77	& 	0	&	53	\\
\instance{Belotti80} (30)	&	11	&	96	&	52	&	11	&	96	&	0	&	52	&	11	&	96	&	52	&	11	&	96	& 	0	&	52	\\
\midrule																											
\instance{Boland80} &	16	&	507	&	46	&	15	&	520	&	13	&	46	&	18	&	507	&	46	&	16	&	507	& 	27	&	46	\\
&	9	&	267	&	23	&	6	&	268	&	8	&	37	&	10	&	267	&	23	&	9	&	267	& 	24	&	23	\\
&	26	&	668	&	17	&	21	&	689	&	11	&	17	&	26	&	668	&	17	&	27	&	668	& 	17	&	17	\\
&	16	&	531	&	19	&	11	&	415	&	13	&	19	&	17	&	531	&	19	&	17	&	531	& 	1	&	19	\\
&	14	&	465	&	22	&	11	&	400	&	12	&	18	&	14	&	465	&	22	&	13	&	465	& 	2	&	22	\\
\cmidrule(r){2-15}	
&\bf 		16	&\bf 		488	&\bf 		25	&\bf 		13	&\bf 		458	&\bf 	11	&\bf 		27	&\bf 		17	&\bf 		488	&\bf 		25	&\bf 		17	&\bf 		488	&\bf 	14	&\bf 		25	\\
\cmidrule(r){2-15}	
\instance{Boland160}	&	430	&	3133	&	13	&	387	&	2944	&	20&	13	&	444	&	3133	&	13	&	445	&	3133& 	85	&	13	\\
&	564	&	2543	&	12	&	483	&	2437	&	20	&	12	&	549	&	2543	&	12	&	544	&	2543	& 	76	&	12	\\
&	241	&	1781	&	13	&	276	&	2303	&	17	&	20	&	233	&	1781	&	13	&	239	&	1781& 	129	&	13	\\
&	782	&	3646	&	15	&	814	&	3768	&	20	&	15	&	763	&	3646	&	15	&	777	&	3646	& 	82	&	15	\\
&	302	&	2021	&	17	&	291	&	2086&	25	&	13	&	291	&	2021	&	17	&	301	&	2021	& 	107	&	17	\\
\cmidrule(r){2-15}	
&\bf 		464	&\bf 		2625	&\bf 		14	&\bf 		450	&\bf 		2708	&\bf 	20	&\bf 		15	&\bf 		456	&\bf 		2625	&\bf 		14	&\bf 		461	&\bf 		2625&\bf 	96	&\bf 		14	\\
\cmidrule(r){2-15}	
\instance{Boland320}	&	13019	&	10862	&	10	&	16403	&	17004	&	39	&	63	&	13009	&	10862	&	10	&	13355	&	10862	& 	390	&	10	\\
&	22572	&	15924	&	8	&	22102	&	17575	&	32	&	8	&	22931	&	15924	&	8	&	22306	&	15924	& 	409	&	8	\\
&	22006	&	14403	&	9	&	24181	&	21072	&	40	&	75	&	21820	&	14403	&	9	&	22153	&	14403	& 	79	&	9	\\
&	21831	&	16990	&	10	&	22486	&	18319	&	35	&	12	&	21837	&	16990	&	10	&	20380	&	16990	& 	357	&	10	\\
&	15981	&	13597	&	9	&	13840	&	12569	&	36	&	9	&	15277	&	13597	&	9	&	14204	&	13597	& 	157	&	9	\\
\cmidrule(r){2-15}	
&\bf 		19082	&\bf 		14355	&\bf 		9	&\bf 		19802	&\bf 		17308	&\bf 	36	&\bf 		33	&\bf 		18975	&\bf 		14355	&\bf 		9	&\bf 		18480	&\bf 		14355	&\bf 	278	&\bf 		9	\\
\cmidrule(r){2-15}	
\instance{Boland16} &	2	&	32	&	5	&	1	&	32	&	1	&	5	&	2	&	32	&	5	&	2	&	32	& 	0	&	5	\\
&	3	&	49	&	11	&	2	&	47	&	1	&	11	&	2	&	49	&	11	&	2	&	49	& 	0	&	11	\\
&	7	&	125	&	27	&	5	&	123	&	1	&	27	&	7	&	125	&	27	&	6	&	125	& 	0	&	27	\\
&	10	&	183	&	25	&	8	&	183	&	1	&	25	&	10	&	183	&	25	&	10	&	183	& 	0	&	25	\\
\cmidrule(r){2-15}	
&\bf 		5	&\bf 		97	&\bf 		17	&\bf 		4	&\bf 		96	&\bf 	1	&\bf 		17	&\bf 		5	&\bf 		97	&\bf 		17	&\bf 		5	&\bf 		97	&\bf 	0	&\bf 		17	\\
\cmidrule(r){2-15}	
\instance{Boland25} &	14	&	162	&	14	&	13	&	183	&	1	&	14	&	13	&	162	&	14	&	14	&	162	& 	0	&	14	\\
&	25	&	283	&	15	&	22	&	289	&	1	&	15	&	26	&	283	&	15	&	25	&	283	& 	0	&	15	\\
&	40	&	429	&	13	&	33	&	422	&	1	&	13	&	39	&	429	&	13	&	40	&	429	& 	0	&	13	\\
&	43	&	437	&	20	&	41	&	466	&	1	&	20	&	44	&	437	&	20	&	43	&	437	& 	0	&	20	\\
\cmidrule(r){2-15}	
&\bf 		31	&\bf 		328	&\bf 		16	&\bf 		27	&\bf 		340	&\bf 	1	&\bf 		16	&\bf 		31	&\bf 		328	&\bf 		16	&\bf 		31	&\bf 		328	&\bf 	0	&\bf 		16	\\
\cmidrule(r){2-15}	
\instance{Boland50} &	395	&	1343	&	14	&	341	&	1409	&	1	&	14	&	397	&	1343	&	14	&	397	&	1343& 	0	&	14	\\
&	754	&	1952	&	17	&	606	&	1890	&	1	&	17	&	766	&	1952	&	17	&	772	&	1952& 	0	&	17	\\
&	1427	&	2593	&	9	&	1249	&	2437	&	1	&	9	&	1382	&	2593	&	9	&	1357	&	2593	& 	0	&	9	\\
&	1740	&	3386	&	15	&	615	&	1622	&	1	&	15	&	1754	&	3386	&	15	&	1702	&	3386	& 	0	&	15	\\
\cmidrule(r){2-15}	
&\bf 	1079	&\bf 		2319	&\bf 		14	&\bf 		703	&\bf 		1840	&\bf 	1	&\bf 		14	&\bf 		1074	&\bf 		2319	&\bf 		14	&\bf 		1057	&\bf 		2319	&\bf 	0	&\bf 		14	\\
\cmidrule(r){2-15}	
\bottomrule
\end{tabular}
\end{table}

Notice from Table \ref{table_presolve} that the results for duality fixing show the opposite pattern for the \instance{Boland320} instances than for all other instances. This is due to the fact that, for an unknown reason, fixing several variables during presolve had a negative impact on preprocessing, causing many fewer solutions to be discovered during this phase and therefore having an overall negative impact on the rest of the BB procedure. We felt though that the positive impact duality fixing had on the other instances sets warranted leaving this feature on for the remainder of our tests. Also observe from  Table \ref{table_presolve} that the exploitation of neither singleton nor dominating columns had any significant impact on the overall BB procedure. We found that this was mainly due to the fact that there were very few occurrences of either of these types of columns. \revision{Singleton columns did not fix any variables values, and so there is no column ``Fixed'' for it in the table}. We opted to turn off the exploitation of singleton columns for the remainder of our tests, but we left on the exploitation of dominating columns. Our reasoning here was that singleton columns have no impact on BB that extends beyond presolve, while dominating columns result in disjunctions from which we can generate global cutting planes. Hence, we left on the exploitation of dominating columns in order to test the impact of generating these cuts in later tests. 

\subsection{Preprocessing}\label{preprocess_test}

In our next test we examined the impact of the two preprocessing techniques discussed in Section \ref{sec:presolve}, as well as a hybrid method we derived as a combination of the two presented procedures. In our initial implementation of this test we used each of these methods with $\rho$ assigned each integer value in $[0,5]$. Recognize from Algorithms \ref{alg5} and \ref{alg6} that each of the proposed preprocessing procedures are designed so that the total number of Pareto solutions computed should have a positive correlation with the value of $\rho$. We determined that \textsc{ProprocesingMethod1} performed poorly for $\rho \leq 1$, \textsc{ProprocesingMethod2} performed poorly for $\rho \geq 2$ and the hybrid method performed poorly in general. Hence, we do not report results for these procedures. We also discovered that the impact of $\rho$ on overall solution time varied with the size of the instance solved. As a result, we also implemented modified preprocessing procedures in which the value of $\rho$ is automatically computed as a function of the size of an instance. 
For each family of instance, the average CPU required to complete BB after employing each of the aforementioned preprocessing strategies is reported in Table \ref{table_preprocessing}. We note that in Table \ref{table_preprocessing} $\rho = v$ indicates that $\rho$ was automatically computed as a function of instance size.

\begin{table}[h!]
\small
\centering
\caption{Experiment 2 -- Impact of preprocessing techniques. Solution time (sec.) to optimality}\label{table_preprocessing}
\begin{tabular}{l|rrrrr|rrr}
\toprule
\multirow{2}{*}{Instance} & \multicolumn{5}{c|}{\textsc{PreprocessingMethod1}} & \multicolumn{3}{c}{\textsc{PreprocessingMethod2}}\\
 & $\rho = 2$ & 3 & 4 & 5 & $v$ & $\rho = 0$ & 1 & $v$\\
\midrule
\instance{Belotti60} (30)	&	4	&	4	&	8	&	8	&	4	&	4	&	5	&	5\\
\instance{Belotti80} (30)	&	11	&	11	&	18	&	18	&	11	&	11	&	12	&	12\\
\midrule																							
\instance{Boland80} &	10	&	11	&	10	&	10	&	11	&	15	&	16	&	15	\\
&	5	&	7	&	6	&	7	&	8	&	6	&	7	&	7	\\
&	20	&	18	&	17	&	18	&	22	&	21	&	26	&	25	\\
&	16	&	17	&	16	&	17	&	16	&	11	&	12	&	13	\\
&	7	&	6	&	7	&	7	&	11	&	11	&	12	&	12	\\
\cmidrule(r){2-9}
&\bf 		12	&\bf 		12	&\bf 		11	&\bf 		12	&\bf 		13	&\bf 		13	&\bf 		14	&\bf 		14\\
\cmidrule(r){2-9}
\instance{Boland160} &	388	&	299	&	218	&	219	&	298	&	401	&	383	&	393	\\
&	300	&	265	&	266	&	263	&	335	&	487	&	516	&	506	\\
&	177	&	144	&	125	&	125	&	158	&	282	&	291	&	287	\\
&	549	&	552	&	541	&	557	&	548	&	816	&	880	&	862	\\
&	171	&	185	&	158	&	158	&	443	&	302	&	280	&	290	\\
\cmidrule(r){2-9}
&\bf 		317	&\bf 		289	&\bf 		262	&\bf 		264	&\bf 		356	&\bf 		458	&\bf 		470	&\bf 		467\\
\cmidrule(r){2-9}
\instance{Boland320} &	11036	&	8619	&	6398	&	6232	&	9561	&	16480	&	16544	&	16636	\\
&	15099	&	16278	&	16210	&	16142	&	14963	&	22319	&	21181	&	21246	\\
&	9433	&	10421	&	9615	&	9840	&	10675	&	24151	&	21878	&	21788	\\
&	14379	&	16642	&	16446	&	16427	&	16253	&	22837	&	24763	&	24384	\\
&	10303	&	10706	&	10779	&	10811	&	10602	&	14422	&	14440	&	14449	\\	
\cmidrule(r){2-9}
&\bf 		12050	&\bf 		12533	&\bf 		11890	&\bf 		11891	&\bf 		12411	&\bf 		20042	&\bf 		19761	&\bf 		19701\\
\cmidrule(r){2-9}
\instance{Boland16} &	1	&	2	&	2	&	3	&	1	&	2	&	3	&	3	\\
&	2	&	2	&	3	&	3	&	2	&	2	&	3	&	3	\\
&	5	&	5	&	6	&	6	&	5	&	5	&	7	&	7	\\
&	8	&	9	&	10	&	9	&	8	&	8	&	10	&	9	\\
\cmidrule(r){2-9}
&\bf 		4	&\bf 		4	&\bf 		5	&\bf 		5	&\bf 		4	&\bf 		4	&\bf 		6	&\bf 		5\\
\cmidrule(r){2-9}
\instance{Boland25} &	11	&	10	&	10	&	10	&	12	&	14	&	16	&	17	\\
&	18	&	19	&	18	&	18	&	18	&	23	&	26	&	25	\\
&	22	&	22	&	22	&	23	&	24	&	31	&	37	&	37	\\
&	50	&	49	&	50	&	53	&	52	&	40	&	44	&	43	\\
\cmidrule(r){2-9}
&\bf 		25	&\bf 		25	&\bf 		25	&\bf 		26	&\bf 		27	&\bf 		27	&\bf 		31	&\bf 		31\\
\cmidrule(r){2-9}
\instance{Boland50} &	278	&	293	&	196	&	198	&	335	&	342	&	354	&	356	\\
&	633	&	546	&	456	&	464	&	663	&	583	&	689	&	678	\\
&	990	&	1110	&	743	&	708	&	945	&	1250	&	1848	&	1852	\\
&	599	&	2217	&	1325	&	1382	&	1325	&	625	&	2054	&	2001	\\
\cmidrule(r){2-9}
&\bf 	625	&\bf 		1042	&\bf 		680	&\bf 		688	&\bf 		817	&\bf 		700	&\bf 		1236	&\bf 		1222\\
\cmidrule(r){2-9}
\bottomrule
\end{tabular}
\end{table}

Observe from Table \ref{table_preprocessing} that 
although variants of \textsc{PreprocessingMethod2} performed well for smaller instances, the same is not true for larger instances. \textsc{PreprocessingMethod1}, on the other hand, performed quite well on all instances. Notice, however, that values of $\rho$ near two performed quite well for small instances while values near five performed extremely poorly. On the other hand, for larger instances values of $\rho$ near five seem to outperform almost every other procedure. Due to the consistent performance of the variant of \textsc{PreprocessingMethod1} with $\rho = 2$, we opted to use this approach for the remainder of our tests.

\subsection{Probing and Pareto Branching}

The next test we performed was designed to examine the utility of the variable probing procedure used directly after preprocessing and at each node prior to branching, and the Pareto branching that we perform when $\OS$ fathoming results in disjoint feasible regions of $\OS$. The results of this experiment are given in Table \ref{table_probing}.

\begin{table}
\small
\centering
\caption{Experiment 3 -- Impact of Probing and Pareto branching}\label{table_probing}
\begin{tabular}{l | rrr | rrr | rr | rr}
\toprule
\multirow{2}{*}{Instance} & \multicolumn{3}{c|}{All Off} & \multicolumn{3}{c|}{Initial Probing} & \multicolumn{2}{c|}{Probing in Branching} & \multicolumn{2}{c}{Pareto Branching}\\
& Time & Nodes & $\bar{\mathcal{G}}_{0}$ & Time & Nodes & $\bar{\mathcal{G}}_{0}$ & Time & Nodes & Time & Nodes  \\ 
\midrule
\instance{Belotti60} (30)	&	4	&	72	&	48	&	4	&	74	&	49	&	3	&	48	&	4	&	72	\\
\instance{Belotti80} (30)	&	11	&	98	&	49	&	11	&	97	&	49	&	8	&	62	&	10	&	94	\\
\midrule																							
\instance{Boland80} &	10	&	368	&	12	&	10	&	366	&	12	&	5	&	160	&	11	&	388	\\
&	5	&	256	&	35	&	6	&	254	&	35	&	5	&	210	&	7	&	260	\\
&	20	&	647	&	19	&	20	&	647	&	19	&	7	&	283	&	17	&	621	\\
&	14	&	598	&	45	&	16	&	581	&	45	&	9	&	393	&	15	&	526	\\
&	7	&	302	&	18	&	7	&	284	&	18	&	5	&	157	&	7	&	332	\\
\cmidrule(r){2-11}
&\bf 		11	&\bf 		434	&\bf 		26	&\bf 		12	&\bf 		426	&\bf 		26	&\bf 		6	&\bf 		241	&\bf 		11	&\bf 		425	\\
\cmidrule(r){2-11}
\instance{Boland160} &	393	&	3185	&	19	&	349	&	2815	&	19	&	125	&	1082	&	259	&	2394	\\
&	309	&	1713	&	20	&	338	&	1743	&	20	&	122	&	625	&	290	&	1948	\\
&	171	&	1466	&	5	&	178	&	1433	&	5	&	91	&	651	&	149	&	1551	\\
&	547	&	2982	&	8	&	547	&	3016	&	8	&	201	&	1249	&	488	&	3595	\\
&	167	&	1196	&	28	&	168	&	1154	&	28	&	78	&	570	&	153	&	1447	\\
\cmidrule(r){2-11}
	&\bf 		318	&\bf 		2108	&\bf 		16	&\bf 		316	&\bf 		2032	&\bf 		16	&\bf 		123	&\bf 		835	&\bf 		268	&\bf 		2187	\\
\cmidrule(r){2-11}
\instance{Boland320} &	10951	&	10391	&	6	&	11061	&	10673	&	6	&	3099	&	3882	&	7120	&	8292	\\
&	14601	&	12827	&	6	&	15038	&	12954	&	6	&	5012	&	5329	&	11358	&	12004	\\
&	9402	&	7626	&	12	&	9316	&	7598	&	12	&	3173	&	3380	&	7571	&	8072	\\
&	14065	&	12161	&	6	&	14542	&	12528	&	6	&	5583	&	5679	&	11685	&	13181	\\
&	9991	&	9900	&	5	&	9850	&	9930	&	5	&	2664	&	3462	&	6555	&	8380	\\
\cmidrule(r){2-11}
	&\bf 		11802	&\bf 		10581	&\bf 		7	&\bf 		11962	&\bf 		10737	&\bf 		7	&\bf 		3906	&\bf 		4346	&\bf 		8858	&\bf 		9986	\\
\cmidrule(r){2-11}
\instance{Boland16} &	1	&	29	&	5	&	1	&	28	&	5	&	1	&	28	&	1	&	47	\\
&	2	&	54	&	12	&	2	&	56	&	12	&	1	&	43	&	2	&	63	\\
&	5	&	128	&	42	&	5	&	124	&	42	&	3	&	104	&	6	&	163	\\
&	7	&	165	&	12	&	7	&	168	&	12	&	5	&	129	&	9	&	199	\\
\cmidrule(r){2-11}
&\bf 		4	&\bf 		94	&\bf 		18	&\bf 		4	&\bf 		94	&\bf 		18	&\bf 		3	&\bf 		76	&\bf 		4	&\bf 		118	\\
\cmidrule(r){2-11}
\instance{Boland25} &	11	&	157	&	32	&	11	&	159	&	32	&	7	&	130	&	10	&	175	\\
&	18	&	343	&	36	&	18	&	337	&	36	&	12	&	259	&	23	&	445	\\
&	23	&	370	&	64	&	29	&	505	&	64	&	15	&	284	&	26	&	379	\\
&	50	&	764	&	76	&	52	&	765	&	76	&	33	&	545	&	38	&	580	\\
\cmidrule(r){2-11}
&\bf 	25	&\bf 	409	&\bf 	52	&\bf 	28	&\bf 	442	&\bf 	52	&\bf 		17	&\bf 		305	&\bf 		24	&\bf 		395	\\
\cmidrule(r){2-10}
\instance{Boland50} &	278	&	1501	&	33	&	304	&	1660	&	33	&	165	&	1063	&	292	&	1831	\\
&	614	&	2318	&	44	&	749	&	2799	&	44	&	499	&	1862	&	585	&	2857	\\
&	948	&	2966	&	22	&	1101	&	3367	&	22	&	600	&	2188	&	704	&	2949	\\
&	559	&	2083	&	60	&	2038	&	5349	&	60	&	1001	&	2583	&	438	&	2135	\\
\cmidrule(r){2-11}
&\bf 	600	&\bf 	2217	&\bf 	40	&\bf 	1048	&\bf 	3294	&\bf 	40	&\bf 	566	&\bf 	1924	&\bf 	505	&\bf 	2443	\\
\cmidrule(r){2-11}
\bottomrule
\end{tabular}
\end{table}

Observe from Table \ref{table_probing} that when utilizing probing directly after preprocessing, in many cases the total CPU time and number of nodes processed increased. Surprisingly, however, performing the same probing procedure prior to branching at each node had an extremely positive impact on the overall performance of BB, significantly lowering total CPU time and the number of explored nodes. We also found that Pareto branching had an overall positive impact on BB performance. For the remainder of our tests we opted to cease probing directly after preprocessing, but to still employ probing during branching as well as Pareto branching. 

\subsection{Exploiting $\OS$ Gaps and Comparing with Triangle Splitting}

We now present the results of an experiment designed to test the performance of our BB procedure against that of the triangle splitting (TS) method of \citep{boland2015acriterion}. For this experiment we solved all the same instances we used in our previous tests and employed two variants of our BB procedure, one in which we utilized the $\OS$ splitting procedure we discussed in Section \ref{sec:osgap} and one in which we utilized our standard implementation. The results of this test are given in Table \ref{table_gaps_ts}. Our standard BB procedure outperformed the triangle splitting method on all but one set of instances, while our $\OS$ splitting procedure outperformed the triangle splitting method on all sets of instances except one. Also recognize that the total CPU times associated with our $\OS$ splitting procedure are always comparable with those of our standard procedure. We point out that there were many more substantial gaps between solutions to exploit after preprocessing for the Belotti* instances than for the Boland* instances. This is the reason that there is a drastic reduction in total number of nodes processed when using $\OS$ splitting on the \instance{Belotti*} instances but not the \instance{Boland*} instances. \revision{We also did a parallel implementation of the $\OS$ splitting procedure and observed some reduction in the CPU times, which suggests that parallelising this procedure can further improve the BB algorithm.} 

\begin{table}[h!]
\small
\centering
\caption{Experiment 4 -- Impact of $\OS$ Gap Splitting and comparison with Triangle Splitting (TS)}\label{table_gaps_ts}
\begin{tabular}{l| rr| rr| r}
\toprule
\multirow{2}{*}{Instance} & \multicolumn{2}{c|}{Standard BB} & \multicolumn{2}{c|}{BB with $\OS$ Gaps} & TS\\
& Time & Nodes & Time & Nodes & Time  \\ 
\midrule
\instance{Belotti60}	(30)	&	3	&	49	&	4	&	33	&	9	\\
\instance{Belotti80}	(30)	&	7	&	64	&	8	&	44	&	20	\\
\midrule	
\instance{Boland80}	&	6	&	205	&	5		&	205	&	44	\\
&	3	&	203	&	3	&	138	&	29	\\
&	7	&	326	&	7	&	261	&	46	\\
&	7	&	300	&	6	&	262	&	48	\\
&	3	&	165	&	3	&	165	&	32	\\
\cmidrule(r){2-6}
&\bf	 5	&\bf	 240	&\bf	 5	&\bf	 206	&\bf 	40	\\
\cmidrule(r){2-6}
\instance{Boland160} &	79	&	914	&	85	&	886	&	320	\\
&	97	&	734	&	105	&	749	&	335	\\
&	67	&	692	&	60	&	668	&	267	\\
&	180	&	1631	&	188	&	1626	&	677	\\
&	56	&	691	&	48	&	573	&	258	\\
\cmidrule(r){2-6}	
&\bf 	96	&\bf 	932	&\bf 	97	&\bf 	900	&\bf 	371	\\
\cmidrule(r){2-6}
\instance{Boland320}	&	2048	&	3391	&	2055	&	3371	&	3800	\\
&	3213	&	4568	&	3333	&	4743	&	6219	\\
&	1981	&	3135	&	1957	&	3164	&	5035	\\
&	3239	&	5429	&	3328	&	5385	&	5421	\\
&	1755	&	3461	&	1912	&	3697	&	4293	\\
\cmidrule(r){2-6}
&\bf 	2447	&\bf 	3997	&\bf 	2517	&\bf 	4072	&\bf 	4954	\\
\cmidrule(r){2-6}
\instance{Boland16} &	1	&	39	&	1	&	39	&	4	\\
&	2	&	47	&	1	&	47	&	5	\\
&	2	&	94	&	4	&	128	&	10	\\
&	5	&	133	&	5	&	133	&	13	\\
\cmidrule(r){2-6}
&\bf 	3	&\bf 	78	&\bf 	3	&\bf 	87	&\bf 	8	\\
\cmidrule(r){2-6}
\instance{Boland25} &	6	&	137	&	6	&	119	&	19	\\
&	14	&	325	&	9	&	215	&	30	\\
&	13	&	258	&	16	&	347	&	39	\\
&	22	&	397	&	24	&	433	&	51	\\
\cmidrule(r){2-6}
&\bf 	14	&\bf	 279	&\bf 	14	&\bf 	279	&\bf 	35	\\
\cmidrule(r){2-6}
\instance{Boland50} &	158	&	1156	&	137	&	961	&	159	\\
&	374	&	2058	&	306	&	1754	&	262	\\
&	484	&	2240	&	371	&	1795	&	346	\\
&	990	&	3843	&	977	&	3369	&	475	\\
\cmidrule(r){2-6}
&\bf 	502	&\bf 	2324	&\bf 	448	&\bf 	1970	&\bf 	311	\\
\cmidrule(r){2-6}
\bottomrule
\end{tabular}
\end{table}

\subsection{Approximations of the Pareto Set}\label{sec:approximation_computation}

\citet{boland2015acriterion} measured the time it takes the Triangle Splitting method to compute an approximate Pareto set having the property that the hypervolume gap between valid primal and dual bounds implied by this approximate set is less than 2\%. We repeat this experiment for our BB procedure, though we note that the primal and dual bounds we utilize are significantly different than those used in \citep{boland2015acriterion}. We measure this gap directly after the completion of our preprocessing procedure, and then each time 25 nodes are processed during BB. We cease the procedure if: (i) BB terminates with the true Pareto set, or (ii) the hypervolume gap is less than 2\%. In this experiment we also report Hausdorff gap measurements, as described in Section \ref{sec:gap}. Additionally, for comparison we include certain results as reported in \citep{boland2015acriterion}.

\begin{table}[h!]
\footnotesize
\centering
\caption{Experiment 5 -- Obtaining approximate Pareto sets}\label{table_appx}
\begin{tabular}{l|rr|rrrrrr|r}
\toprule
 & \multicolumn{8}{c|}{Standard BB} & TS\\
\multirow{3}{*}{Instance}  & \multicolumn{2}{c|}{Preprocessing} & \multicolumn{6}{c|}{Until $\mathcal{HV}_{s^*} \le 2\%$} & \\
& $\mathcal{HV}_0$ & $\bar{\mathcal{G}}_{0}$ & Time & \% Time & Nodes & \% Nodes & $\mathcal{HV}_{s^*}$ & $\bar{\mathcal{G}}_{s^*}$ &  \% Time \\ 
\midrule
\instance{Belotti60} (30)	&	21.7	&	62.4	&	3	&	100	&	49	&	98	&	0.2	&	4.9	&	--	\\
\instance{Belotti80} (30)	&	25.7	&	66.8	&	7	&	100	&	62	&	98	&	0.2	&	3.3	&	--	\\
\midrule																					
\instance{Boland80}	&	1.6	&	12.3	&	1	&	18	&	0	&	0	&	1.6	&	12.3	&	12	\\
&	3.5	&	34.7	&	3	&	74	&	100	&	49	&	2.0	&	16.3	&	9	\\
&	2.9	&	19.0	&	2	&	31	&	25	&	8	&	1.8	&	10.9	&	4	\\
&	49.0	&	67.1	&	5	&	74	&	225	&	75	&	1.2	&	10.7	&	6	\\
&	2.2	&	18.2	&	1	&	47	&	25	&	15	&	1.2	&	11.2	&	7	\\
\cmidrule(r){2-10}																					
& \bf	11.8	& \bf	30.3	& \bf	2	& \bf	49	& \bf	75	& \bf	29	& \bf	1.6	& \bf	12.3	& \bf	7.6	\\
\cmidrule(r){2-10}																					
\instance{Boland160}	&	2.6	&	18.6	&	16	&	21	&	75	&	8	&	1.2	&	8.8	&	2.30	\\
&	1.9	&	20.0	&	8	&	9	&	0	&	0	&	1.9	&	20.0	&	3.85	\\
&	1.2	&	4.7	&	4	&	6	&	0	&	0	&	1.2	&	4.7	&	1.50	\\
&	0.8	&	7.9	&	10	&	5	&	0	&	0	&	0.8	&	7.9	&	0.61	\\
&	9.2	&	28.4	&	20	&	35	&	150	&	22	&	1.8	&	8.4	&	2.90	\\
\cmidrule(r){2-10}																					
& \bf	3.1	& \bf	15.9	& \bf	12	& \bf	15	& \bf	45	& \bf	6	& \bf	1.4	& \bf	10.0	& \bf	2.23	\\
\cmidrule(r){2-10}																					
\instance{Boland320}	&	1.1	&	6.4	&	52	&	3	&	0	&	0	&	1.1	&	6.4	&	0.21	\\
&	0.5	&	5.9	&	82	&	3	&	0	&	0	&	0.5	&	5.9	&	0.23	\\
&	0.5	&	12.3	&	78	&	4	&	0	&	0	&	0.5	&	12.3	&	0.26	\\
&	0.5	&	5.9	&	80	&	2	&	0	&	0	&	0.5	&	5.9	&	0.23	\\
&	0.4	&	5.5	&	72	&	4	&	0	&	0	&	0.4	&	5.5	&	0.22	\\
\cmidrule(r){2-10}																					
& \bf	0.6	& \bf	7.2	& \bf	73	& \bf	3	& \bf	0	& \bf	0	& \bf	0.6	& \bf	7.2	& \bf	0.23	\\
\cmidrule(r){2-10}																					
\instance{Boland16} &	0.6	&	5.3	&	1	&	68	&	0	&	0	&	0.6	&	5.3	&	--	\\
&	1.2	&	11.9	&	0	&	26	&	0	&	0	&	1.2	&	11.9	&	--	\\
&	3.1	&	42.3	&	2	&	74	&	25	&	27	&	1.0	&	18.6	&	--	\\
&	2.2	&	12.3	&	2	&	34	&	25	&	19	&	1.6	&	11.4	&	--	\\
\cmidrule(r){2-10}																					
& \bf	1.7	& \bf	18.0	& \bf	1	& \bf	50	& \bf	13	& \bf	11	& \bf	1.1	& \bf	11.8	&	--	\\
\cmidrule(r){2-10}																					
\instance{Boland25}	&	4.0	&	32.2	&	5	&	82	&	75	&	55	&	1.2	&	9.7	&	--	\\
&	67.3	&	79.4	&	13	&	94	&	175	&	54	&	1.8	&	19.9	&	--	\\
&	83.0	&	87.8	&	7	&	55	&	75	&	29	&	1.5	&	28.5	&	--	\\
&	91.2	&	93.7	&	11	&	50	&	100	&	25	&	1.9	&	22.0	&	--	\\
\cmidrule(r){2-10}																					
& \bf	61.4	& \bf	73.3	& \bf	9	& \bf	70	& \bf	106	& \bf	41	& \bf	1.6	& \bf	20.0	&	--	\\
\cmidrule(r){2-10}																					
\instance{Boland50}	&	2.6	&	33.3	&	15	&	10	&	25	&	2	&	1.9	&	28.9	&	--	\\
&	3.8	&	44.3	&	152	&	41	&	325	&	16	&	2.0	&	28.5	&	--	\\
&	1.9	&	21.8	&	6	&	1	&	0	&	0	&	1.9	&	21.8	&	--	\\
&	26.6	&	75.1	&	191	&	19	&	275	&	7	&	1.7	&	24.0	&	--	\\
\cmidrule(r){2-10}																					
& \bf	8.7	& \bf	43.6	& \bf	91	& \bf	18	& \bf	156	& \bf	6	& \bf	1.9	& \bf	25.8	&	--	\\
\cmidrule(r){2-10}																					
\bottomrule
\end{tabular}
\end{table}

The results of this experiment are displayed in Table \ref{table_appx} from which we make several observations. For the majority of the Boland* instances, the hypervolume gap is already less than 2\% after preprocessing, before BB even begins. This is evidence that these instances are relatively easy. Recall Figure \ref{pareto_sets}, and notice that for the \instance{Boland80} instance the boundary of the dual bound at the root node is very close to the Pareto set. This is further evidence of the ease of these instances. In contrast to this, notice from Table \ref{table_appx} that for the Belotti* instances, it takes over 75\% of the total BB time in order to obtain a hypervolume gap of less than 2\%. We note that Table \ref{table_appx} also shows that the triangle splitting method is able to determine an approximate solution with a hypervolume gap of less than 2\% in less time, relative to the total solution time. 

\renewcommand{\dagger}{\ast}

\subsection{MIPLIB Instances}\label{sec:MIPLIB}
\revision{
Due to the successful results we obtained using our BB procedure on instances from the literature, we designed our second set of tests to measure the performance of our algorithm on a more challenging set of instances. For this we utilized a set of single objective MILP instances available from the MIPLIB 2017 library \cite[p]{miplib2017pre}. We chose only instances that were feasible, mixed-integer, contained at most 1000 total decision variables, and were marked easy. For each instance, we generated two secondary objective functions as follows:

\begin{enumerate}
\item[(r)] For each $i\in \{1,\dots,m+n\}$ the coefficient $c^2_i$ is randomly generated using the uniform distribution over the closed interval $\left[-\left\vert\max_i c^1_i\right\vert, \left\vert\max_i c^1_i\right\vert\right]$.
\item[(n)] We set $c^2_i = -c^1_i$.
\end{enumerate}

After generation of these instances we did some preliminary testing and discarded instances for which: (i) the Pareto set was a singleton, or (ii) the second objective was unbounded, or (iii) the MILP associated with either $f_1$ or $f_2$ took over 8 hours to solve. In the end, 104 instances remained for final testing (2 each, originating from 52 single objective MILP instances).

A primary reason for generating this additional set of instances is the relative ease with which single objective MILPs were solved throughout the solution process, both during the execution of triangle splitting and our \bb{}, when using the previously considered test sets. As such, in our first analysis of these new instances we set a variety of node limits on single objective MILPs solved during our \bb{} (other than the two specified on line 4 of Algorithm \ref{alg7}, of course). By limiting the number of nodes processed during each single objective MILP solve, we hoped to increase the speed of the overall \bb{} procedure while still being able to exploit useful dual bound information at each node. For initial tests, we set node limits of 10, $10^2$, $10^3$, $10^4$, and $\infty$ and compared the overall solution time for \bb{} on all 104 instances, with a maximum execution time of 8 hours. Surprisingly, the best performing node limits were $10^4$ and $\infty$. Hence, we opted to leave the single objective MILP node limit off for the remainder of our analysis. We did note, however, that on some instances, single objective MILPs took significant time to solve even when relatively few nodes were explored in order to do so. Thus, we opted to solve each instance again, this time with an overall time limit imposed when solving each single objective MILP. For this test we utilized time limits of 15, 30, 45, 60, 300, 1800, and $\infty$ seconds. In this case, the limits that appeared to produce the best results were 30 and 300 seconds, with 300 seconds having a slight advantage. We therefore imposed a single objective MILP time limit of 300 seconds when conducting our final round of tests.
}

\begin{table}
\centering
\caption{Experiment 6 -- Solution time (sec.) or duality gap $\mathcal{G}_{s^{*}}$ (\%)  for biobjective instances generated from MIPLIB 2017}\label{table_results}
\tiny
\begin{tabular}{lrrr || lrrr}
\toprule
& \multicolumn{2}{c}{Branch-and-bound} & Triangle Splitting && \multicolumn{2}{c}{Branch-and-bound} & Triangle Splitting \\
\cmidrule(r){2-3}\cmidrule(r){6-7}
Instance & Standard & Gap Splitting && Instance & Standard & Gap Splitting  \\
\midrule
\instance{22433}		&	 5.53	&	 3.93	&	  5.36	& \instance{neos-1425699}		&	$\dagger$	&	\errortwo	&	$\dagger$	\\
			&	21.31	&	59.43	&	139.40 &&	$\dagger$	&	$\dagger$	&	$\dagger$	\\[8pt]
\instance{23588}		&	60.69	&	50.60	&	 91.47	& \instance{neos-1430701}		&	5797.07	&	27292.20	&	2086.53	\\
			&	24616.54	&	25092.16	&	790.16 &&	1.87	&	1.90	&	2.64	\\[8pt]
\instance{assign1-5-8}		&	25599.20	&	25104.12	&	1867.39	& \instance{neos-1442119}		&	\percentGap{0.40}	&	\percentGap{0.38}	&	22552.6	\\
		&	2035.07	&	2031.20	&	1415.81 &&	326.48	&	328.68	&	1534.18		\\[8pt]
\instance{b-ball}		&	0.60	&	0.31	&	0.06	& \instance{neos17}		&	$\dagger$	&	$\dagger$	&	25302.70	\\
		&	0.01	&	0.01	&	0.02 &&	$\dagger$	&	$\dagger$	&	$\dagger$	\\[8pt]
\instance{beavma}		&	4178.27	&	3615.45	&	$\dagger$	& \instance{neos-3610041-iscar}		&	731.45	&	142.51	&	8.14 \\
		&	$\dagger$	&	$\dagger$	&	$\dagger$ &&	3.37	&	3.42	&	3.91	\\[8pt]
\instance{blend2}		&	5232.48	&	5274.62	&	$\dagger$ &	\instance{neos-3610051-istra}		&	2002.65	&	880.83	&	37.81\\
		&	$\dagger$	&	$\dagger$	&	28405.80 &&	$\dagger$	&	\errortwo	&	6.34	\\[8pt]
\instance{ci-s4}		&	\percentGap{2.85}	&	\percentGap{2.86}	&	$\dagger$	& \instance{neos-3610173-itata}		&	12251.52	&	\errortwo	&	75.74 \\
		&	13511.01	&	15171.00	&	$\dagger$	& &	\percentGap{0.71}	&	\percentGap{9.98}	&	25.90 \\[8pt]
\instance{dcmulti}		&	725.25	&	728.54	&	173.88	& \instance{neos-3611447-jijia}		&	2934.49	&	1354.50	&	37.68\\
		&	$\dagger$	&	$\dagger$	&	0.76	& &	$\dagger$	&	\percentGap{1.17}	&	19.07 \\[8pt]
\instance{exp-1-500-5-5}		&	$\dagger$	&	$\dagger$	&	$\dagger$ & \instance{neos-3611689-kaihu}		&	398.39	&	\errortwo	&	53.70	\\
		&	5.59	&	5.56	&	21.36 &&	$\dagger$	&	$\dagger$	&	18.41	\\[8pt]
\instance{fastxgemm-n2r6s0t2}		&	\percentGap{11.37}	&	\percentGap{11.94}	&	178.48 & \instance{neos5}		&	285.47	&	286.68	&	265.46	\\
		&	1419.83	&	1414.09	&	0.6 & &	2613.94	&	3046.41	&	273.85	\\[8pt]
\instance{flugpl}		&	0.48	&	0.44	&	1.46 & \instance{Neos-5192052-neckar}		&	0.79	&	0.79	&	$\dagger$	\\
		&	0.02	&	0.02	&	0.04 & &	9.72	&	9.71	&	0.06	\\[8pt]
\instance{gen}		&	72.61	&	74.01	&	15.57 & \instance{Nexp-50-20-1-1}		&	$\dagger$	&	\errortwo	&	1318.04	\\
		&	\percentGap{24.58}	&	\percentGap{35.37}	&	0.32 &&	\errorone &	\errorone	&	28.08	\\[8pt]
\instance{gr4x6}		&	0.69	&	0.56	&	1.12	& \instance{noswot}		&	407.97	&	\errortwo 	&	261.77 \\
		&	$\dagger$	&	$\dagger$	&	5.29 & &	272.75	&	195.12	&	171.7	\\[8pt]
\instance{ic97\_potential}		&	\percentGap{13.05}	&	\percentGap{12.16}	&	$\dagger$ & \instance{nsa}		&	\percentGap{14.47}	&	\errortwo	&	$\dagger$	\\
		&	21411.78	&	21597.36	&	$\dagger$ &&	\percentGap{0.94}	&	\percentGap{8.06}	&	$\dagger$	\\[8pt]
\instance{ic97\_tension}		&	\percentGap{4.52}	&	\percentGap{4.55}	&	$\dagger$ & \instance{opt1217}		&	0.21	&	0.21	&	0.15	\\
		&	9.10	&	9.12	&	53.19 &&	0.10	&	0.09	&	0.02	\\[8pt]
\instance{k16x240b}		&	$\dagger$	&	$\dagger$	&	$\dagger$	& \instance{prod1}		&	\percentGap{3.71}	&	19240.07	&	978.34\\
		&	390.86	&	664.73	&	$\dagger$	&&	20.29	&	23.34	&	27.52\\[8pt]
\instance{markshare\_4\_0}		&	3665.90	&	112.33	&	9841.90 & \instance{prod2}		&	\percentGap{3.04}	&	\percentGap{1.98}	&	25214.90	\\
		&	\errorone	&	\errorone	&	768.88 &&	283.95	&	283.36	&	183.98	\\[8pt]
\instance{markshare\_5\_0}		&	\percentGap{0.93}	&	$\dagger$	&	27098.70 & \instance{qiu}		&	863.70	&	1002.68	&	18267.50	\\
		&	\percentGap{0.27}	&	\percentGap{0.27}	&	6767.01 &&	$\dagger$	&	\percentGap{7.80}&	$\dagger$	\\[8pt]
\instance{mas74}		&	27118.70	&	20170.72	&	19182.20	& \instance{r50x360}		&	$\dagger$	&	$\dagger$	&	$\dagger$\\
		&	486.88	&	489.42	&	$\dagger$	& &	\errorone	&	\errorone	&	$\dagger$\\[8pt]
\instance{mas76}		&	$\dagger$	&	\percentGap{50}	&	$\dagger$	& \instance{ran12x21}		&	$\dagger$	&	\percentGap{4.68}	&	$\dagger$\\
		&	$\dagger$	&	$\dagger$	&	$\dagger$ & &	$\dagger$	&	$\dagger$	&	$\dagger$	\\[8pt]
\instance{mik\_250\_20\_75\_1}		&	$\dagger$	&	\percentGap{2.90}	&	$\dagger$ & \instance{ran13x13}		&	8732.94	&	8420.79	&	$\dagger$	\\
		&	$\dagger$	&	\percentGap{7.79}	&	20.57 & &	$\dagger$	&	$\dagger$	&	$\dagger$	\\[8pt]
\instance{mik\_250\_20\_75\_2}		&	$\dagger$	&	$\dagger$	&	$\dagger$ & \instance{ran14x18-disj-8}		&	\percentGap{1.70}	&	\percentGap{1.75}	&	$\dagger$	\\
		&	$\dagger$	&	$\dagger$	&	14.62 & &	$\dagger$	&	\percentGap{< 0.01}	&	$\dagger$	\\[8pt]
\instance{mik\_250\_20\_75\_3}		&	$\dagger$	&	$\dagger$	&	$\dagger$ & \instance{rout}		&	1191.12	&	1095.00	&	2795.34	\\
		&	$\dagger$	&	$\dagger$	&	15.60 & &	$\dagger$	&	$\dagger$	&	$\dagger$	\\[8pt]
\instance{mik\_250\_20\_75\_4}		&	$\dagger$	&	\percentGap{1.18}	&	$\dagger$ & \instance{sp150x300d}		&	$\dagger$	&	\percentGap{0.48}	&	$\dagger$	\\
		&	\percentGap{1.01}	&	\percentGap{1.42}	&	139.28 & &	$\dagger$	&	\percentGap{2.93}	&	39.27	\\[8pt]
\instance{mik\_250\_20\_75\_5}		&	$\dagger$	&	$\dagger$	&	$\dagger$ & \instance{timtab1}		&	\percentGap{15.19}	&	$\dagger$	&	$\dagger$	\\
		&	$\dagger$	&	\percentGap{6.38}	&	17.69 & &	< \percentGap{0.01}	&	\percentGap{2.41}	&	$\dagger$	\\[8pt]
\instance{misc07}		&	$\dagger$	&	$\dagger$	&	3816.02 & \instance{timtab1CUTS}		&	\percentGap{23.55}	&	\percentGap{48.75}	&	$\dagger$	\\
		&	\percentGap{4.42}	&	\percentGap{16.89}	&	614.38 & &	$\dagger$	&	$\dagger$	&	$\dagger$	\\[8pt]
\bottomrule
\end{tabular}
\end{table}

\begin{table}
\centering
\caption{Number of primal solutions (singletons plus line segments) at termination for MIPLIB test set}
\label{table_size_prim}
\scriptsize
\begin{tabular}{lr || lr || lr || lr}
\toprule

Instance & Sol & Instance & Sol & Instance & Sol & Instance & Sol  \\
\midrule
\instance{	22433	}	&	3	&	\instance{	ic97\_potential	}	&	158	&	\instance{	Neos-1425699	}	&	2415	&	\instance{	nsa	}	&	22	\\
			&	46	&				&	6	&				&	148983	&				&	118	\\
\instance{	23588	}	&	10	&	\instance{	ic97\_tension	}	&	117	&	\instance{	Neos-1430701	}	&	33	&	\instance{	opt1217	}	&	2	\\
			&	234	&				&	4	&				&	2	&				&	2	\\
\instance{	Assign1-5-8	}	&	11	&	\instance{	k16x240b	}	&	2521	&	\instance{	Neos-1442119	}	&	33	&	\instance{	prod1	}	&	38	\\
			&	2	&				&	38	&				&	2	&				&	2	\\
\instance{	B-ball	}	&	4	&	\instance{	Markshare\_4\_0	}	&	126	&	\instance{	neos17	}	&	400	&	\instance{	prod2	}	&	54	\\
			&	2	&				&	2811	&				&	96	&				&	2	\\
\instance{	beavma	}	&	473	&	\instance{	Markshare\_5\_0	}	&	285	&	\instance{	neos-3610041-iscar	}	&	53	&	\instance{	qiu	}	&	600	\\
			&	6235	&				&	4992	&				&	64	&				&	1948	\\
\instance{	blend2	}	&	255	&	\instance{	mas74	}	&	20	&	\instance{	neos-3610051-istra	}	&	52	&	\instance{	r50x360	}	&	1760	\\
			&	1536	&				&	7	&				&	52	&				&	219	\\
\instance{	Ci-s4	}	&	365	&	\instance{	mas76	}	&	230	&	\instance{	neos-3610173-itata	}	&	141	&	\instance{	ran12x21	}	&	3740	\\
			&	4	&				&	89200	&				&	142	&				&	1705	\\
\instance{	dcmulti	}	&	243	&	\instance{	Mik\_250\_20\_75\_1	}	&	2638	&	\instance{	neos-3611447-jijia	}	&	157	&	\instance{	ran13x13	}	&	2586	\\
			&	5	&				&	32	&				&	204	&				&	2923	\\
\instance{	Exp-1-500-5-5	}	&	419	&	\instance{	Mik\_250\_20\_75\_2	}	&	3638	&	\instance{	neos-3611689-kaihu	}	&	78	&	\instance{	ran14x18-disj-8	}	&	2729	\\
			&	4	&				&	32	&				&	204	&				&	2150	\\
\instance{	Fastxgemm-n2r6s0t2	}	&	38	&	\instance{	Mik\_250\_20\_75\_3	}	&	2898	&	\instance{	neos5	}	&	30	&	\instance{	rout	}	&	2	\\
			&	2	&				&	32	&				&	8	&				&	32817	\\
\instance{	flugpl	}	&	87	&	\instance{	Mik\_250\_20\_75\_4	}	&	2584	&	\instance{	Neos-5192052-neckar	}	&	77	&	\instance{	sp150x300d	}	&	103	\\
			&	4	&				&	32	&				&	5	&				&	232	\\
\instance{	gen	}	&	193	&	\instance{	Mik\_250\_20\_75\_5	}	&	2715	&	\instance{	Nexp-50-20-1-1	}	&	118	&	\instance{	timtab1	}	&	324	\\
			&	4	&				&	32	&				&	217	&				&	148	\\
\instance{	gr4x6	}	&	40	&	\instance{	misc07	}	&	629	&	\instance{	noswot	}	&	7	&	\instance{	timtab1CUTS	}	&	238	\\
			&	50	&				&	630	&				&	37	&				&	112	\\
\bottomrule
\end{tabular}
\end{table}

\begin{figure}[h!]
\begin{subfigure}[h]{.55\textwidth}
\centering
\includegraphics[width=.9\textwidth]{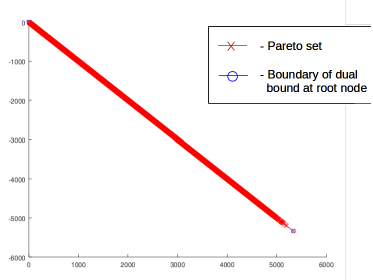}
\caption{Variant (n) of the ``markshare\_5\_0'' instance.}
\label{mark5_v2}
\end{subfigure}%
\begin{subfigure}[h]{.45\textwidth}
\centering
\includegraphics[width=1\textwidth]{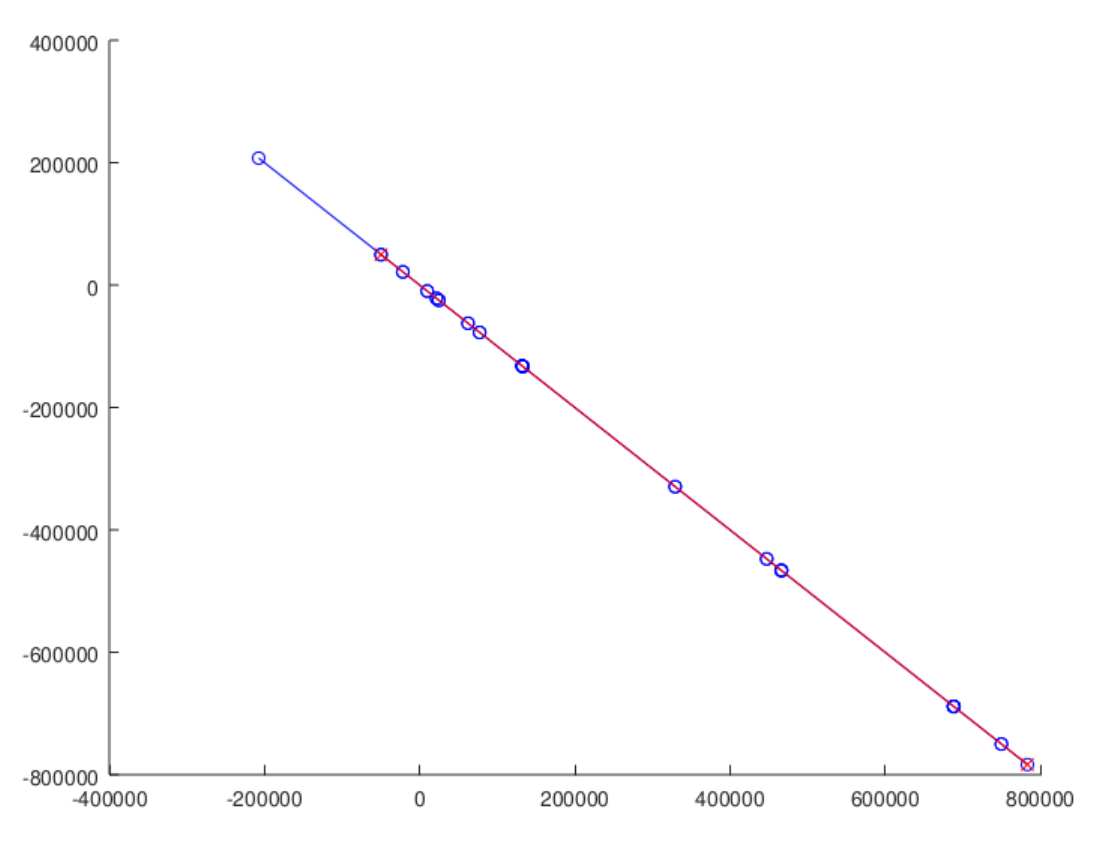}
\caption{Variant (n) of the ``mik\_25\_20\_75\_1'' instance.}
\label{mik1_v2}
\end{subfigure}
\caption{Pareto set and boundary of $\L_0$ for two instances from MIPlib 2017.}
\label{pareto_sets2}
\end{figure}

\revision{
Table~\ref{table_results} gives the results of this experiment, where the two lines for each instance correspond to the (r) and (n) methods, respectively, for generating the second objective function. Of the 104 instances considered, 46 were solved in under 8 hours by the original \bb{} implementation, 44 by the $\OS$ splitting \bb{} variant, and 64 by the triangle splitting method. Additionally, there were 8 instances which were solved in under 8 hours by at least one version of \bb{}, but not by the triangle splitting method, and 23 instances solved in under 8 hours by the triangle splitting method, but not by a \bb{} procedure. When \bb{} failed to complete in 8 hours, we allowed an additional hour for the computation of $\bar{\mathcal{G}}_{s^*}$. This computation was sometimes unsuccessful due to the extreme number of open nodes present upon termination. When it was successful, it is reported in Table~\ref{table_results} using the \percentGap{} symbol, otherwise we display a $\dagger$ symbol. Table~\ref{table_size_prim} gives the number of nondominated solutions stored when our BB terminated either by solving to optimality or reached its time limit.

In all, the results display comparable performance between the BB approaches and the triangle splitting method, though for instances in which there was a relative difference in performance, it was generally large. These discrepancies in performance seem to stem from overall structure of the Pareto set in $\OS$. In particular, triangle splitting appears have superior performance on instances for which either: (i) the total number of Pareto solutions is small, or (ii) most Pareto solutions are supported,  particularly if all Pareto solutions lie along a single line segment in $\OS$. The latter of these two properties was observed on the second variant of the ``markshare\_5\_0'' instance, for example. Additionally, both of these properties were observed in the Pareto sets of the second variant of instances beginning with ``mik,'' for example. To aid in the visualization of this property, we include Figure \ref{pareto_sets2} which displays the Pareto set and boundary of $\L_0$ for the (n) variants of the ``markshare\_5\_0'' and ``mik\_25\_20\_75\_1'' instances. On the other hand, \bb{} appears to have superior performance on instances for which either: (i) the total number of Pareto solutions is large, or (ii) a relatively large percentage of Pareto solutions are unsupported. We also note that occasionally numerical issues caused early termination of \bb{} when solving instances for which all Pareto solutions fall on a single line segment in $\OS$. In particular, for some such instances, the cutting plane generated on line 21 of Algorithm \ref{alg1} fell along this line segment in $\OS$ when the root node was processed. Thus, CPLEX determined that the duality gap here was 0 and terminated \bb{} execution prematurely. This phenomenon was observed on the (n) variant of the ``markshare\_4\_0,'' ``Nexp-50-20-1-1,'' and ``r50x360'' instances. This is indicated in Table \ref{table_results} using the \errorone{} symbol. In addition, the implementation of \bb{} which exploits $\OS$ gaps failed due to numerical issues on a small number of instances. This is indicated in Table \ref{table_results} using the \errortwo{} symbol.
}

\section{Concluding Remarks} \label{sec:conclude}
In this paper, we have introduced a new \bb{} method for solving BOMILP with general integers. For each component of single objective \bb{}, we presented procedure(s) for extending this component to the  biobjective setting. We have also conducted numerous computational experiments. The first several experiments provide insight into the usefulness of each of the algorithms we proposed. The final few experiments compare the performance of our BB procedure and the triangle splitting method \citep{boland2015acriterion}. Our BB procedure outperforms the triangle splitting method on  instances from literature, and performs comparably on large, challenging instances that were developed in this paper. 

Most of the algorithms proposed by us have, in theory, straightforward generalizations to the multiobjective case (MOMILPs). However, having an implementable correct \bb{} for MOMILPs is far from a trivial extension of this work. We point out some important questions that need to be answered in this regard.

\subsection{Extension to multiobjective MILP}
Correct node fathoming is what makes a \bb{} algorithm a correct and exact method. Fathoming by bound dominance is how fathoming mostly occurs in \bb{}. For BOMILP, the bound sets are two-dimensional polyhedra. This greatly simplifies checking bound dominance for BOMILPs since given two line segments, or piecewise linear curves in general, in $\R^{2}$, one can easily identify the dominated portion through pairwise comparisons. The data structure \citep{treestructure} stores nondomimated line segments and efficiently checks if a new line segment is dominated by what is currently stored. This enabled the node processing step in this paper to perform fathoming efficiently. Bound sets for MOMILP are higher-dimensional polyhedra and hence one will require an even more sophisticated data structure to store these sets. Since the local dual bound set at each node is a polyhedron and the global primal bound is a finite union of polyhedra, checking dominance requires checking containment of polyhedra, whose complexity depends on their respective representations, 
and also computing the set difference between the primal and dual bound sets. The set resulting from this set difference would be nonconvex, in general, which begs the question: is there a straightforward way to represent this nonconvex set as a union of polyhedra whose relative interiors are disjoint?  Thus, fathoming and storing nondominated regions for a MOMILP is even more nontrivial. Once these obstacles are overcome, the \bb{} proposed in this paper should extend to a implementable \bb{} for MOMILPs.

{
\newrefcontext[sorting=nyt]
\printbibliography
}

\end{document}